\documentclass[a4paper,10pt,twoside]{article}
\usepackage{authblk}

\usepackage[ngerman,english]{babel}
\usepackage[latin1]{inputenc}
\usepackage{csquotes}
\usepackage{amsfonts,amsmath,amsthm}
\usepackage{empheq}

\usepackage[titletoc,title]{appendix}
\usepackage[backend=bibtex8,doi=false,eprint=false,firstinits=true,isbn=false,style=numeric-comp,url=false,maxnames=99]{biblatex}
\makeatletter
\def\blx@maxline{77}
\makeatother
\bibliography{bibl.bib}
\DeclareRedundantLanguages{english,german,french}{english,german,ngerman,french}

\usepackage{cases}
\usepackage{mathabx}
\usepackage{mathtools}
\usepackage{xfrac}
\usepackage{enumerate}
\usepackage{fancyhdr}
\usepackage[usenames,dvipsnames,svgnames,table]{xcolor}
\usepackage[colorlinks]{hyperref}
\definecolor{blue75}{rgb}{0,0,.75}
\definecolor{green75}{rgb}{0,.75,0}
\hypersetup{colorlinks=true, urlcolor=blue75,linkcolor=blue75,citecolor=green75,pdfstartview=FitB,bookmarksopen=true,bookmarksopenlevel=1}

\usepackage[a4paper, left=2.5cm, right=2.5cm, top=2.5cm,bottom=2cm]{geometry}
\usepackage{constants}
\newcommand{\parenthezises}[1]{\arabic{#1}}
\newconstantfamily{C}{
symbol=C,
format=\parenthezises,
reset={section}
}
\newconstantfamily{M}{
symbol=M,
format=\parenthezises,
reset={section}
}
\usepackage{enumitem}

\allowdisplaybreaks
\begin{document}

\newcommand{\cmg}{\color{magenta}}
\newcommand{\D}{\mathbb{D}}
\newcommand{\R}{\mathbb{R}}
\newcommand{\N}{\mathbb{N}}
\newcommand{\p}{\mathbb{P}}

\def\avint{\mathop{\,\rlap{-}\!\!\int}\nolimits} 
\def\diam{\operatorname{diam}}
\def\dist{\operatorname{dist}}
\def\ess{\operatorname{ess}}
\def\inner{\operatorname{int}}
\def\osc{\operatorname{osc}}
\def\sign{\operatorname{sign}}
\def\supp{\operatorname{supp}}

\newtheorem{Assumptions}{Assumptions}[section]
\newtheorem{Theorem}[Assumptions]{Theorem}
\newtheorem{Corollary}[Assumptions]{Corollary}
\newtheorem{Convention}[Assumptions]{Convention}
\newtheorem{Definition}[Assumptions]{Definition}
\newtheorem{Example}[Assumptions]{Example}
\newtheorem{Lemma}[Assumptions]{Lemma}
\newtheorem{Notation}[Assumptions]{Notation}
\newtheorem{Remark}[Assumptions]{Remark}

\numberwithin{equation}{section}
\title{The Malliavin derivative and compactness: an application to a degenerate PDE-SDE coupling}
\author{Anna~Zhigun}
\renewcommand\Affilfont{\itshape\small}
\affil{Technische Universität Kaiserslautern, Felix-Klein-Zentrum für Mathematik\\ Paul-Ehrlich-Str. 31, 67663 Kaiserslautern, Germany\\
  e-mail: {zhigun@mathematik.uni-kl.de}}
\date{}
\maketitle
\begin{abstract}
Compactness is one of the most versatile tools in the analysis of nonlinear PDEs and systems. Usually, compactness is established by means of some embedding theorem between functional spaces. Such theorems, in turn, rely on appropriate estimates for a function and its derivatives. While a similar result based on simultaneous estimates for the Malliavin and weak Sobolev derivatives is available for the Wiener-Sobolev spaces, it seems that it has not yet been widely used in the analysis of highly nonlinear parabolic problems with stochasticity. In the present work we apply this result in  order to study compactness,  existence  of global solutions, and, as a by-product, the convergence of a semi-discretisation scheme for a prototypical degenerate PDE-SDE coupling.
\end{abstract}
\section{Introduction}
Compactness is one of the most versatile tools in the analysis of nonlinear equations and systems. Usually, it is established by means of a compactness criterion for a particular functional space. Well-known examples  include: the Rellich-Kondrachov theorem, the Lions-Aubin lemma,  the Arzel\`{a}-Ascoli and Fr\'echet-Kolmogorov theorems.  Such results rely on appropriate estimates of a function and its classical or weak derivatives or, more generally, of its increments. These theorems are most helpful instruments in the study of  deterministic differential equations and systems in H\"older and Sobolev spaces.
A closely related result was established for  Wiener-Sobolev spaces  by Bally and Saussereau \cite{BallySaussereau}.  It is based on simultaneous estimates for the Malliavin and weak Sobolev derivatives. It seems, however, that this criterion has not yet been widely used in the analysis of highly nonlinear parabolic problems with stochasticity. Indeed, standard approaches mostly rely on some  monotonicity of the elliptic part (see \cite{KryRoz,LiuRoeckner,PrRoeckner} and references therein) which often  fails to hold for strongly coupled systems. The approach based on a priori estimates and a compact embedding has several advantages. For instance, it allows to treat rather general classes of complicated problems by approximating them with better-studied, more regular ones, following the so-called compactness method \cite{Lions}. In the present work we adopt this scheme. Namely, we apply the result by Bally and Saussereau in  order to study compactness and  existence  of global solutions for a  prototypical degenerate PDE-SDE coupling. The proof of existence is based on the semi-discretisation method. As a by-product, it justifies the convergence of a semi-discretisation scheme for our problem.

This paper is organised as follows. First,   we introduce our model system  in {\it Section  \ref{model}}, fix some notations in {\it Section \ref{not}},  and state the main results in {\it Section \ref{problem}}. 
We then establish in {\it Sections \ref{aprioric}-\ref{aprioriy}} a set of uniform a priori estimates for the solutions of our system. While estimates in {\it Subsections \ref{aprioricbasic}-\ref{EstDrDif}} and {\it Section \ref{aprioriy}} are rather standard, those in {\it Subsections \ref{secPMTr}-\ref{SecPsi}} are new and more involved.  The a priori  estimates lie at the core of the proof of compactness in {\it Section \ref{SecComp}}. In {\it Section \ref{SecSemi}} we introduce a spatial discretisation scheme and study its convergence in the nondegenerate case. Finally, we use compactness  in order to prove the existence of global solutions to the original degenerate problem  in {\it Section \ref{existence}}.

\textbf{Acknowledgment.} The author expresses her thanks to Wolfgang Bock and Christina Surulescu (Technische Universität Kaiserslautern) for  stimulating discussions. {The idea to consider a model in the form of a PDE-SDE coupling 
was proposed by Christina Surulescu in the context of \cite{KlSS}.}
\section{The model}\label{model}
Let $T$ be a positive number and $O$ be a smooth bounded domain in $\R^N$. Also, let $\left(\Omega,{\mathcal F},({\mathcal F}_t\right)_{t\in[0,T]},\p)$ be a filtered probability space  on which is defined a Wiener process $(W(t))_{t\in[0,T]}$. We assume that $\left({\mathcal F}_t\right)_{t\in[0,T]}$ is the usual completion of the natural filtration of $(W(t))_{t\in[0,T]}$. 
In this setting, we consider the following system of a random porous medium equation (random PME, RPME)
\begin{subequations}\label{PM}
 \begin{empheq}[left={\text{in } \Omega\ \empheqlbrace\,}]{align}
&\partial_t \beta(c)=\Delta c+f(c,y)&&\text{ in }(0,T]\times O,\label{c}\\
 &c=0\qquad (\partial_{\nu}c=0)&&\text{ in }(0,T]\times \partial O,\label{bc}\\
 &c=c_0&&\text{ in }\{0\}\times O\label{inic0}
\end{empheq}
\end{subequations}
and an It\^{o} SDE
\begin{subequations}\label{SDE}
 \begin{empheq}[left={\text{in } O\ \empheqlbrace\,}]{align}
 &dy=a(y)\,d W+b(c,y)\,dt &&\text{ in } (0,T]\times\Omega,\label{y}\\
 &y=y_0&&\text{ in } \{0\}\times\Omega.\label{iniy0}
\end{empheq}
\end{subequations}
Couplings of a PDE, however with a linear diffusion, with an SDE have lately emerged in the multiscale tumor modelling  \cite{KlSS}. In that work, the PDE- and  SDE-variables represented, respectively, the intra- and extracellular proton dynamics in a tumor. In the present case, the variables $c$ and $y$ could be seen, e.g., as the tumor density and concentration of the intracellular protons, respectively (see also \cite{Hiremath2015176,HirSur2016} for models based on PDE-RODE-couplings). Thereby, the  (one-dimensional) Wiener process in SDE  \eqref{y} captures some  stochastic fluctuations in the intracellular proton dynamics (see \cite{KlSS} and references therein). Through the coupling terms $f$ and $b$, these fluctuations influence indirectly  the dynamics of the cancer cells on the macroscale as well.

The class of functions $\beta$ considered in the present study (see {\it Assumptions \ref{Assump1}} below) includes as a particular case
\begin{align}
\beta(c)=c^{\frac{1}{m}}\text{ for some }m>1.\nonumber
\end{align}
For such $\beta$, we can transform the  macroscopic PDE \eqref{c} into the equation 
\begin{align}
 \partial_t u=\Delta u^m+f(u^m,y),\nonumber
\end{align}
where 
\begin{align*}
  u:=c^{\frac{1}{m}}.
\end{align*}
Thus, switching to the new variable $u$, we regain a RPME with a source term, written in the conventional form. 
Porous medium equations are standard examples of degenerate-diffusion equations. As in our recently proposed deterministic models, \cite{ZSU,ZSH},  degeneracy  accounts for a finite speed of  propagation of a tumor. While in \cite{ZSU,ZSH} we assumed the diffusion coefficient to be degenerate not only in $c$, but also in another variable, here we consider a simpler case, with the diffusion being of the porous medium type.
\section{Basic notation and functional spaces}\label{not}
We denote $\R^+:=(0,\infty)$, $\R^+_0:=[0,\infty)$. 

For a Lebesgue measurable set $E$ we denote by $|E|$ its Lebesgue measure. The space dimension depends on the context. The integral average of an integrable function $f:E\rightarrow\R$ is defined via
\begin{align*}
 \avint_{E}f(x)\,dx:=\frac{1}{|E|}\int_{E}f(x)\,dx.
\end{align*}

We assume a smooth bounded domain $O\subset\R^N$, $N\in\N$, to be given. The  outward unit normal vector on the boundary of $O$ we denote by $\nu$.

The derivative of a function $u$ of one real variable is denoted by $u'$. Partial derivatives in the  classical or distributional sense with respect to a variable  $z$ are denoted by $\partial_{z}$. The variable $z$ can, for example, be  the 'time' variable $t\in \R^+_0$ or a component of the 'spatial' variable $x\in \overline{O}$.  Further, $\nabla$  and $\Delta$ stand for the spatial gradient and Laplace operator, respectively.

We assume the reader to be familiar with the standard $L^p$, Sobolev, and H\"older spaces and their usual properties, as well as with the  more general $L^p$ spaces of functions with values in general Banach spaces, and with anisotropic spaces. In particular, for relatively open w.r.t. $(0,T]\times O$  (not necessarily cylindrical) sets $Q$, we need the spaces
\begin{align}
 &W^{(1,2),q}(Q):=\left\{u=u(t,x)\in L^q(Q):\ \partial_t u, \partial_{x_i}u,\partial_{x_i}\partial_{x_j}u\in L^q(Q)\text{ for }i,j=1,\dots,N\right\}\nonumber
 \end{align}
 with a norm defined via
 \begin{align}
 &\|u\|_{W^{(1,2),q}(Q)}:=\left(\|u\|_{L^q(Q)}^q+\|\partial_t u\|_{L^q(Q)}^q+\sum_{i=1}^{N}\left\|\partial_{x_i}u\right\|_{L^q(Q)}^q+\sum_{i,j=1}^{N}\left\|\partial_{x_i}\partial_{x_j}u\right\|_{L^q(Q)}^q\right)^{\frac{1}{q}}.\nonumber
\end{align}
We recall that the H\"older coefficient for a H\"older exponent $\gamma\in(0,1)$ and a real-valued function $w$ defined in a set $A\subset \R^k$, $k\in\N$, is given by 
\begin{align*}
 |w|_{C^{\gamma}(A)}:=\underset{x,y\in A,\ x\neq y}{\sup}\frac{|w(x)-w(y)|}{|x-y|^{\gamma}}.
\end{align*}

By $\partial A$ we denote the topological boundary of a set $A\subset \R^k$, $k\in\N$. For a set $Q\subset[0,T]\times \overline{O}$, we call 
\begin{align}
 \partial Q\backslash (\{T\}\times O)\nonumber
\end{align}
the parabolic boundary of $Q$.

Further, let a filtered probability space $\left(\Omega,{\mathcal F},({\mathcal F}_t\right)_{t\geq0},\p)$ on which is defined a Wiener process $(W(t))_{t\geq0}$ be given.  We assume that $\left({\mathcal F}_t\right)_{t\geq0}$ is the usual completion of the natural filtration of $(W(t))_{t\geq0}$.   The corresponding It\^o differential is denoted by $dW$. We presuppose that the reader is familiar with some It\^{o} and Malliavin calculi. In particular, we assume such standard results as: the It\^{o} isometry (see, e.g., \cite[Chapter 1 Theorem 7.1]{Ma}) for $p=2$), a version of the Burkholder-Gundy-Davis inequality as stated in  \cite[Chapter 1 Theorem 7.1]{Ma}, as well as the Kolmogorov's continuity criterion (see, e.g., \cite[Chapter 1 Theorem 3.1]{Ch} and the subsequent remark) to be known. 
For $T>0$, the Malliavin derivative of an ${\mathcal F}_T$-measurable $c$ is denoted by $(D_r c)_{r\in[0,T]}$. 
We make use of the following properties of the Malliavin derivative: the chain rule \cite[Proposition 1.2.3]{nualart2006malliavin} and a result on the weak differentiability of solutions to    It\^{o} SDEs (see, e.g., \cite[Theorem 2.2.1]{nualart2006malliavin} and the subsequent observation). We refer to \cite{Ma,Ch,nualart2006malliavin} for more details on the calculus for stochastic processes.

To ease the notation while dealing with purely PDE (SDE) properties which hold $\p$-a.s. in $\Omega$ (a.e. in $O$), 
we sometimes drop the dependence upon variable $\omega$ (variable $x$) and write, for example, $c$ instead of $c(\cdot,\omega,\cdot)$ ($c(\cdot,\cdot,x)$). 
Moreover, for a stochastic process $u:[0,T]\times\Omega\rightarrow V$, where $V$ is a space of functions defined for  $x\in\overline{D}$,
we often write $u(t)$ instead of $u(t,\cdot)$. 

Finally, we make the following two useful conventions. Firstly, for all indices $i$,  $C_i$ or $\alpha_i$ denotes a non-negative constant or, alternatively, a non-negative function, 
which is non-decreasing in each of its arguments. Secondly, we assume that  the reappearing numbers $T,r,t,h_1$, and $h_2$ always satisfy
\begin{align}
 0<r<t\leq T,\qquad 0<h_1\leq T-t,\qquad 0<h_2\leq t-r.\nonumber
\end{align}

\section{Problem setting and main result}\label{problem}
We make the following assumptions on the problem parameters. 
\begin{Assumptions}\label{Assump1}~
\begin{enumerate}
 \item\label{conO} $O$ is a smooth bounded domain in $\R^N$, $N\in \N$.
 \item\label{conbet} Function $\beta:\R^+_0\rightarrow\R^+_0$ satisfies for some constants $m_1>1$, $m_2\geq 1$, and $M,\mu>0$ the conditions
\begin{alignat}{3}
 &\beta\text{ has an inverse function  }\beta^{(-1)},&&\qquad\beta(0)=0,\label{betam1}\\
 &\beta\in C^1(\R^+)\cap C(\R^+_0),&&\qquad\beta'>0\text{ in }\R^+,\label{betaderp}\\
 &\beta'\text{ is decreasing in }\R^+,&&\label{betadecr}\\
 &(\beta')^{-1}\in C^{1-\frac{1}{m_1}}(\R^+_0),&&\label{betamm}\\
 &c^{1-\frac{1}{m_2}}\beta'(c)\leq \mu\quad \text{ for }c\in(0,M].&&\label{betamm2}
\end{alignat}
\item\label{confab} Functions $f,b:\R^+_0\times\R^+_0\rightarrow\R$, and $a:\R^+_0\rightarrow\R$ satisfy the conditions
\begin{alignat}{3}
  &\underset{c,y\in\R^+_0}{\sup}\frac{f(c,y)}{\beta(c)+1}<\infty,&&\qquad\label{fsublin}\\
  &f\in C^1(\R^+_0\times\R^+_0),&&\qquad \partial_c f,\partial_y f\in C(\R^+_0;C_b(\R^+_0)),\label{fLip}\\
  &a\in C^1(\R^+_0),&&\qquad a'\in C_b(\R^+_0),\label{aLip}\\
  &b\in C^1(\R^+_0\times\R^+_0),&&\qquad \partial_c b,\partial_y b\in C(\R^+_0;C_b(\R^+_0)).\label{bLip}
\end{alignat}
Moreover, $f,a$, and $b$ also satisfy
\begin{align}
  &f(0,y)\geq0\text{ for all }y\in\R^+_0,\qquad a(0)=0,\qquad b(c,0)\geq0\text{ for all }c\in\R^+_0.\label{pospar}
\end{align}
\item The initial data $c_0,y_0:\overline{O}\rightarrow\R^+_0$ satisfy 
\begin{align}
&c_0\in L^{\infty}(O)\cap H^1(O),\label{c0}\\
 &y_0\in H^1(O), \label{y0}
\end{align}
as well as 
\begin{align}
  c_0,y_0\geq0\quad\text{ a.e. in }O.\label{inipos}
\end{align}
\end{enumerate}
\end{Assumptions}
\begin{Definition}[Weak-strong solution]\label{Defsol}
Let {\it Assumptions \ref{Assump1}} be satisfied. Let $T>0$. 
We call a pair  $(c,y):[0,T]\times \overline{O}\times\Omega\rightarrow\R^+_0\times \R^+_0$ a \underline{weak-strong} local solution  of system
 \eqref{PM}-\eqref{SDE} if  it holds that:
\begin{enumerate}
\item\label{Defsol1} $(c,y):[0,T]\times\Omega\rightarrow L^2(O)$ is an adapted process;
\item  $c\in L^{\infty}(\Omega;L^{\infty}((0,T);L^{\infty}(O)))$, $\partial_t c\in L^{\infty}(\Omega;L^2((0,T);L^2(O)))$;
\item $c\in L^{\infty}(\Omega;L^{\infty}((0,T); H^1_0(O)))$\quad ($c\in L^{\infty}(\Omega;L^{\infty}((0,T); H^1(O)))$);
\item $D_r c\in L^2((r,T);L^2(\Omega;L^2(O)))$, $\partial_t D_r\beta(c)\in L^2((r,T)\times\Omega;H^{-2}(O))$ for all $r\in[0,T]$;
\item $y\in L^{\infty}((0,T);L^2(\Omega;H^1(O)))\cap C^{\frac{1}{2}}([0,T];L^q(\Omega;L^q(O)))$ for some $q>2$;
\item\label{Defsol5} $D_r y\in L^{\infty}((r,T);L^2(\Omega;L^2(O)))$ for all $r\in[0,T]$;
\item $(c,y)$ satisfies 
 \begin{subequations}
 \begin{empheq}[left={\empheqlbrace\,}]{align}
 &-\left(\beta\left(c_0\right),v\right)_{L^2(O)}\xi(0)-\int_0^T\left(\beta\left(c\right)(s),v\right)_{L^2(O)}\frac{d\xi}{dt}(s)\,ds\nonumber\\
 &=-\int_0^T\left(\nabla c(s),\nabla  v\right)_{L^2(O)}\xi(s)+\left(f\left(c,y\right)(s),v\right)_{L^2(O)}\xi(s)\,ds\nonumber\\&\text{for all }v\in H_0^1(O)\qquad (v\in H^1(O))\quad\text{ for all } \xi\in C_0^{\infty}[0,T)\text{ a.s. in }\Omega,\label{inteqc}\\
 &y(t)=y_0+\int_0^t a(y)(s)\,d W(s)+\int_0^tb(c,y)(s)\,ds\quad\text{ in }L^2(O)\text{ for all }t\in[0,T]\text{ a.s. in }\Omega.\label{inteqy}
\end{empheq}
\end{subequations}
\end{enumerate}
If $(c,y)$ is a local solution for all $T>0$, then we call it a 
global solution.
\end{Definition}
\begin{Remark}[Weak-strong solution]
 We call a solution $(c,y)$ from {\it Definition \ref{Defsol}} a weak-strong solution since it satisfies the PDE \eqref{PM} in a weak PDE-sense and the SDE \eqref{SDE} in the strong SDE-sense.
\end{Remark}
\begin{Remark}[Continuity of  sample paths]\label{hconti}~
\begin{enumerate}
  \item Conditions $c\in L^{\infty}(\Omega;L^{\infty}((0,T);H^1(O)))$ and $\partial_t c\in L^{\infty}(\Omega;L^2((0,T);L^2(O)))$ imply  that  $c(\cdot,\cdot,\omega)\in H^1((0,T);L^2(O))$ a.s. in $\Omega$. Standard result \cite[Chapter 1 Lemma 1.2]{Lions} yields that $c(\cdot,\cdot,\omega)\in C([0,T];L^2(O))$ a.s. in $\Omega$.
 \item  Since $q>2$, condition $y\in C^{\frac{1}{2}}([0,T];L^q(\Omega;L^q(O)))$ implies that $y(\cdot,\cdot,\omega)\in C([0,T];L^q(O))$ a.s. in $\Omega$. This is a direct consequence of the Kolmogorov's continuity criterion (see \cite[Chapter 1 Theorem 3.1]{Ch} and the subsequent remark). 
 Below we choose $q$ as in the Sobolev embedding theorem, i.e., such that $H^1(O)\subset L^q(O)$. 
\end{enumerate}
\end{Remark}
\noindent For arbitrary $T>0$ and  $R=(R_0,R_1)$, $R_0,R_1>0$, we introduce the parameter set
\begin{alignat}{3}
 &\mathcal{P}(T,R)\nonumber\\
 :=&\Big\{(\beta,f,a,b,c_0,y_0)|\ \text{as in {\it Assumptions \ref{Assump1}} and such that }\|c_0\|_{L^{\infty}(O)},\  \underset{c,y\in\R^+_0}{\sup}\frac{f(c,y)}{\beta(c)+1}\leq R_0,\nonumber\\
 &\ \left\|(\beta')^{-1}\right\|_{C^{1-\frac{1}{m_1}}\left[0,R_2(T,R_0)\right]},\ \underset{c\in(0,R_2(T,R_0)]}{\sup}c^{1-\frac{1}{m_2}}\beta'(c),\
 \left\|\left(\partial_c f,\partial_y f,\partial_c b,\partial_y b\right)\right\|_{C_b\left(\left[0,R_2(T,R_0)\right]\times\R^+_0\right)}  \leq R_1,\nonumber\\
 &\ \ \|a'\|_{C_b(\R^+_0)},\ b(0,0),\ \|(c_0,y_0)\|_{H^1(O)}\leq R_1 \Big\},\nonumber
 \end{alignat}
 where
 \begin{align}
 R_2(T,R_0):=\beta^{(-1)}\left(e^{TR_0}(\beta(R_0)+1)-1\right)\label{R2}
\end{align}
is an upper bound for solution component $c$ if the parameters belong to $\mathcal{P}(T,R)$ (see estimate \eqref{cR2} below).  
We define the corresponding  sets of solutions and of their components via 
\begin{align}
 &{\mathcal{U}}(T,R):=\{(c,y)|\ (c,y)\text{ is a weak-strong solution corresponding to } (\beta,f,a,b,c_0,y_0)\in {\mathcal{P}}(T,R) \},\nonumber\\
 &{\mathcal{C}}(T,R):=\{c|\ (c,y)\in {\mathcal{U}}(T,R)\},\nonumber\\
 &{\mathcal{Y}}(T,R):=\{y|\ (c,y)\in {\mathcal{U}}(T,R)\}.\nonumber
\end{align}
Now we are ready to formulate our compactness result:
\begin{Theorem}[Compactness]\label{CompTh}
Let  {\it Assumptions \ref{Assump1} \ref{conO}.-\ref{confab}.} be satisfied. Then for all $T>0$ and $R=(R_0,R_1)$, $R_0,R_1>0$, it holds that
 \begin{align}
  {\mathcal{C}}(T,R)\text{ and }{\mathcal{Y}}(T,R)\text{ are precompact in } L^2((0,T)\times O\times\Omega).\nonumber
 \end{align}
\end{Theorem}
\noindent We prove this theorem in {\it Section \ref{SecComp}} and then use it in {\it Section \ref{existence}} in order to establish the following result which deals with the existence of solutions:
\begin{Theorem}[Existence of a weak-strong solution]\label{ExistTh}
Let {\it Assumptions \ref{Assump1}} be satisfied. Then there exists a weak-strong global solution in terms of {\it Definition \ref{Defsol}}  to system \eqref{PM}-\eqref{SDE}.
\end{Theorem}
\begin{Remark}[Uniqueness] The uniqueness of solutions to \eqref{PM}-\eqref{SDE} holds as well.  It can be proved in a standard way by exploiting the monotonicity of $\beta$ without requiring the solutions to be differentiable in the Malliavin sense. 
\end{Remark}
\begin{Remark}[Notation]
We make the following useful convention: 
the statement that a quantity (a constant or a function) depends  on the parameters of the problem means that it depends  upon the space dimension $N$,   domain $O$, 
constants $T$, $m_1$, $m_2$, $\mu$, and $q$, and the structure of the initial values $c_0$ and $y_0$ and of the  coefficient functions $\beta,f,a$ and $b$ (the latter means their norms etc. which appear in the definition of the parameter set $\mathcal{P}(T,R)$).

Moreover, dependence upon these parameters is mostly 
{\bf not} indicated in an explicit way.
\end{Remark}
\section{A priori estimates for the RPME \texorpdfstring{\eqref{PM}}{PM}}\label{aprioric}
\subsection{Standard PDE estimates for \texorpdfstring{$c$}{c}}\label{aprioricbasic}
Equation \eqref{c} is a.s. in $\Omega$ a standard PME in $(0,T]\times O$. This allows us to derive in a standard way several basic estimates for $c$ which hold irrespectively of $\omega$. To begin with, we multiply \eqref{c} by $p(\beta(c)+1)^{p-1}$ for an arbitrary $p>1$ and integrate by parts over $O$ using the boundary conditions. We thus obtain with the help of assumptions \eqref{betaderp} and \eqref{fsublin} that 
\begin{align}
 \frac{d}{dt} \|\beta(c)+1\|_{L^p(O)}^p=&-p(p-1)\int_{O}c\beta'(c)(\beta(c)+1)^{p-2}|\nabla c|^2\,dx+p\int_{O}f(c,y)(\beta(c)+1)^{p-1}\,dx\nonumber\\
 \leq&p\Cr{Cf}\|\beta(c)+1\|_{L^p(O)}^p,\label{estp1}
\end{align}
where 
\begin{align}
 \Cl{Cf}:=\underset{c,y\in\R^+_0}{\sup}\frac{f(c,y)}{\beta(c)+1}.\nonumber
\end{align}
Applying the Gronwall lemma to \eqref{estp1} and taking the $\frac{1}{p}$-power on both sides of the resulting inequality, we obtain  that
\begin{align}
 \|\beta(c)(t)+1\|_{L^p(O)}\leq e^{t\Cr{Cf}} \|\beta(c_0)+1\|_{L^p(O)}.\label{estp2}
\end{align}
In the limit as $p\rightarrow\infty$ estimate \eqref{estp2} and the fact that $\beta$ is increasing yield that
\begin{align}
  \|\beta(c)(t)+1\|_{L^{\infty}(O)}\leq& e^{t\Cr{Cf}} \|\beta(c_0)+1\|_{L^{\infty}(O)}\nonumber\\
  \leq&e^{t\Cr{Cf}}\left(\beta\left(\|c_0\|_{L^{\infty}(O)}\right)+1\right).\label{estp3}
\end{align}
Consequently, we arrive with \eqref{estp3} and assumption  \eqref{betam1} at the estimate
\begin{align}
 \|c\|_{L^{\infty}((0,T)\times O\times \Omega)}\leq& \beta^{(-1)}\left(e^{T\Cr{Cf}}\left(\beta\left(\|c_0\|_{L^{\infty}(O)}\right)+1\right)-1\right),\nonumber\\
 \leq&R_2(T,R_0),\label{cR2}
\end{align}
where 
\begin{align}
 R_0:=\max\{\Cr{Cf},\|c_0\|_{L^{\infty}(O)}\},\nonumber
\end{align}
and $R_2$ was defined in \eqref{R2}.
\begin{Remark}
 Due to estimate \eqref{cR2}, it suffices to consider the coefficient functions for $c\in[0,R_2(T,R_0)]$ only. 
\end{Remark}
\noindent
Another standard estimate for $c$ as solution to the PME \eqref{PM} is obtained by multiplying by $\partial_t c$ and integrating by parts over $O$ using the boundary conditions and then over $[0,t]$ for $t\in(0,T]$. This implies due to \eqref{cR2}, assumptions \eqref{betadecr} and \eqref{c0}, and the Young inequality that
\begin{align}
\Cl{C0}\int_0^t\left\|\partial_t c(s)\right\|_{L^2(O)}^2\,ds
 \leq &\beta'\left(\|c\|_{L^{\infty}((0,T)\times O\times \Omega)}\right)\int_0^t\left\|\partial_t c(s)\right\|_{L^2(O)}^2\,ds\nonumber\\
 \leq&\int_0^t\int_{O}\beta'(c)|\partial_t c|^2(s)\,dxds\nonumber\\
 =&-\|\nabla c(t)\|_{(L^2(O))^N}^2+\|\nabla c_0\|_{(L^2(O))^N}^2+\int_0^t\int_{O}f(c,y)\partial_t c(s)\,dx ds\nonumber\\
 \leq&-\|\nabla c(t)\|_{(L^2(O))^N}^2+\frac{\Cr{C0}}{2}\int_0^t\left\|\partial_t c(s)\right\|_{L^2(O)}^2\,ds+\C.\nonumber
\end{align}
Consequently, we obtain that
\begin{align}
 \int_0^t\left\|\partial_t c(s)\right\|_{L^2(O)}^2\,ds,\ \|\nabla c(t)\|_{(L^2(O))^N}^2\leq\C\text{ for all }t\in[0,T].\label{est1}
\end{align}
Altogether, estimates \eqref{cR2} and \eqref{est1} yield the first group of estimates for norms of $c$:
 \begin{empheq}[left={ \text{a.s. in }\Omega\ \empheqlbrace\,}]{align}
 &\left\|c\right\|_{L^{\infty}((0,T)\times O)}\leq \Cl{Cap},\label{cbnd}\\
 &\left\|\nabla c\right\|_{L^{\infty}((0,T);(L^2(O))^N)}\leq \Cr{Cap},\label{nablac}\\
 &\left\|\partial_t c\right\|_{L^2((0,T);L^2(O))}\leq \Cr{Cap}.\label{dtc}
\end{empheq}
In particular, assumption \eqref{fsublin} together with estimate \eqref{cbnd} imply that
\begin{align}
 &\left\|f(c)\right\|_{L^{\infty}((0,T)\times O)}\leq \Cl{Capf}\qquad\text{a.s. in }\Omega.\label{fbnd}
\end{align}

\subsection{Estimate for the Malliavin derivative  \texorpdfstring{$D_rc$}{Drc}}
%
The next step is to apply the Malliavin  derivative operator $D_r$ on both sides of the integrated form of equation \eqref{PM}, the integro-differential equation \eqref{inteqc}. 
Using the chain rule and the locality property ($D_{t_1} F(t_2)=0$ for $t_2<t_1$ for $F$ adapted) of the operator $D_r$, 
we compute that
\begin{align}
 \beta'(c)D_rc(t) =&D_r\beta(c)(t)\nonumber\\
 =&D_r\Delta\int_0^t c(s)\,ds+D_r\int_0^tf(c,y)(s)\,ds\nonumber\\
  =&\Delta\int_r^t D_rc(s)\,ds+\int_r^tD_rf(c,y)(s)\,ds.\label{Dreqc}
 \end{align}
 Multiplying \eqref{Dreqc} by $D_rc(t)$ and integrating by parts over $O$  using the boundary conditions and then over $[r,t]$, we obtain using assumption \eqref{betadecr},  estimate \eqref{cbnd},   and the Young inequality that
 \begin{align}
 &\Cl{C25}\int_r^t\left\|D_rc(s)\right\|_{L^2(O)}^2\,ds\nonumber\\
  \leq&\int_r^t\int_{O}\beta'(c)|D_rc|^2(s)\,dxds\nonumber\\
  =&-\frac{1}{2}\left\|
  \nabla\int_r^t D_rc(s)\,ds\right\|_{(L^2(O))^N}^2+\int_r^t\int_{O}D_rc(s)\int_{r}^sD_rf(c,y)(\tau)\,d\tau dx ds\nonumber\\
  \leq& \frac{\Cr{C25}}{2}\int_r^t\left\|D_rc(s)\right\|_{L^2(O)}^2\,ds+\Cl{C1}\int_r^t\int_{r}^s\left\|D_rf(c,y)(\tau)\right\|_{L^2(O)}^2\,d\tau ds.\nonumber
 \end{align}
 Thus, we arrive at the inequality
\begin{align}
  \int_r^t\left\|D_rc(s)\right\|_{L^2(O)}^2\,ds
  \leq& \C\int_r^t\int_{r}^s\left\|D_rf(c,y)(\tau)\right\|_{L^2(O)}^2\,d\tau ds.\label{dinDrc2}
 \end{align}
 In order to estimate the term on the right hand side of \eqref{dinDrc2}, we use the chain rule and assumption \eqref{fLip}. We thus obtain that
 \begin{align}
  \left\|D_rf(c,y)(t)\right\|_{L^2(O)}^2\leq \Cl{C26}\left\|D_rc(t)\right\|_{L^2(O)}^2+\Cr{C26}\left\|D_ry(t)\right\|_{L^2(O)}^2.\label{fest}
 \end{align}
Combining \eqref{dinDrc2} and \eqref{fest} yields that
\begin{align}
  \int_r^t\left\|D_rc(s)\right\|_{L^2(O)}^2\,ds
  \leq& \Cl{C1_}\int_r^t\int_{r}^s\left\|D_rc(\tau)\right\|_{L^2(O)}^2\,d\tau ds+\Cr{C1_}\int_r^t\int_{r}^s\left\|D_ry(\tau)\right\|_{L^2(O)}^2\,d\tau ds.\label{dinDrc1}
 \end{align} 
Application of  the Gronwall lemma to \eqref{dinDrc1} yields that
\begin{align}
 \int_r^t\left\|D_rc(s)\right\|_{L^2(O)}^2\,ds\leq &\C\int_r^t\int_r^s\left\|D_ry(\tau)\right\|_{L^2(O)}^2\,d\tau\,ds\nonumber\\
 \leq&\Cl{CDr1}\int_r^t\left\|D_ry(s)\right\|_{L^2(O)}^2\,ds\qquad \text{a.s. in }\Omega.\label{estDrc01}
\end{align}
Moreover, estimates  \eqref{fest} and \eqref{estDrc01} imply that
\begin{align}
&\int_r^t\left\|D_rf(c,y)(s)\right\|_{L^2(O)}^2\,ds\leq \Cl{CDrf}\int_r^t\left\|D_ry(s)\right\|_{L^2(O)}^2\,ds\qquad \text{a.s. in }\Omega.\label{estDrf1}
\end{align}
 \subsection{Estimate for  \texorpdfstring{$\partial_t D_r\beta(c)$}{}}
Applying the differential operator $D_r$ to the original equation \eqref{c}, we have that 
\begin{align}
  \partial_t D_r\beta(c)=\Delta D_rc+D_r f(c,y)\label{Dreqc1}
\end{align}
Combining estimates  \eqref{estDrc01} and \eqref{estDrf1} with the Sobolev inequality, we conclude from  \eqref{Dreqc1} that
\begin{align}
 \int_r^t\left\|\partial_t D_r\beta(c)(s)\right\|_{H^{-2}(O)}^2\,ds\leq&\C\int_r^t\left\|D_rc(s)\right\|_{L^2(O)}^2+\left\|D_rf(c,y)(s)\right\|_{H^{-2}(O)}^2\,ds\nonumber\\
 \leq&\C\int_r^t\left\|D_rc(s)\right\|_{L^2(O)}^2+\left\|D_rf(c,y)(s)\right\|_{L^2(O)}^2\,ds\nonumber\\
 \leq&\Cl{CdtDr1}\int_r^t\left\|D_ry(s)\right\|_{L^2(O)}^2\,ds\qquad \text{a.s. in }\Omega.\label{estDrbeta}
\end{align}
 \subsection{Estimate for \texorpdfstring{$D_{r+h_2}c-D_rc$}{}}\label{EstDrDif}
Inserting $r+h_2$ in place of $r$ in \eqref{Dreqc}, we have (recall that $D_{t_1} F(t_2)=0$ for $t_2<t_1$ for $F$ adapted) that
 \begin{align}
  \beta'(c)D_{r+h_2}c(t)
  =&\Delta\int_{r+h_2}^t D_{r+h_2}c(s)\,ds+\int_{r+h_2}^t D_{r+h_2}f(c,y)(s)\,ds,\nonumber\\
  =&\Delta\int_r^t D_{r+h_2}c(s)\,ds+\int_r^t D_{r+h_2}f(c,y)(s)\,ds.\label{DrhDr}
 \end{align}
Subtracting \eqref{DrhDr} and \eqref{Dreqc}, we obtain 
that
 \begin{align}
  \beta'(c)(D_{r+h_2}c-D_r c)(t)
  =&\Delta\int_r^t (D_{r+h_2}c-D_r c)(s)\,ds+\int_r^t (D_{r+h_2}f(c,y)-D_r f(c,y))(s)\,ds.\label{DrhDr2}
 \end{align}
 Arguing for \eqref{DrhDr2} as we just did for \eqref{Dreqc} above, we obtain that
  \begin{align}
 \int_r^t\left\|(D_{r+h_2}c-D_r c)(s)\right\|_{L^2(O)}^2\,ds
 \leq& \Cl{CDr2}\int_r^t\left\|(D_{r+h_2}y-D_r y)(s)\right\|_{L^2(O)}^2\,ds\qquad \text{a.s. in }\Omega.\label{estDrc4}
 \end{align}
\subsection{A regularising transformation for solutions of a PME}\label{secPMTr}
In this sequel, we deal with some purely PDE properties of \eqref{c}. 
It is well-known that  solutions of a PME like \eqref{c}  are, in general, only weak-strong solutions if $c\not\equiv0$ and is not strictly separated from zero. In particular,  $\Delta c$ is generally not even $L^2$-bounded.
Still, it is well-understood \cite{AroBen1979,DiBen83,CaffFri1980,Ivanov86,Ziemer1982} that under reasonable assumptions on $f$ the solution is at least locally H\"older continuous.  
Following our idea from \cite{AZPMReg}, we show how the information on the (local) H\"older continuity of a solution function  can be used in order to transform this function into a smooth one by means of a smooth and strictly increasing function  which depends only upon the parameters of the problem. We believe that this result is of interest by itself.
We then use this transformation for the compactness proof, see {\it Section \ref{SecComp}}. Before we begin with a construction for  solutions of a PME, let us consider a simple motivating example.
\begin{Example}\label{simple} Denote by $B_r$, $r>0$,  the closed ball of radius $r$ in $\R^N$ centred at the origin. Let  
\begin{align*}
 w:B_1\rightarrow[0,1],\qquad w(x)=|x|^{\gamma}\text{ for some }\gamma\in(0,1).
\end{align*}
It is well-known that $w\in C^{2}(B_1\backslash\{0\})\cap C^{\gamma}(B_1)$, but $w\not\in C^2(B_1)$. Since 
\begin{align}
\left|\partial_{x_i}\partial_{x_j}w(x)\right|
\leq \gamma(\gamma-1)|x|^{\gamma-2}
\text{ for all }x\in B_1\backslash\{0\}\nonumber
\end{align}
and 
\begin{align}
\overline{\{w>k\}}=B_{k^{\frac{1}{\gamma}}} \text{ for all }k\in(0,1],\nonumber
\end{align}
we have that
\begin{align}
 \|\partial_{x_i}\partial_{x_j}w\|_{C\left(\overline{\{w>k\}}\right)}\leq\gamma(\gamma-1)k^{1-\frac{2}{\gamma}}
\text{ for all }k\in(0,1].\label{sing}
\end{align}
Set
\begin{align*}
 \varphi:[0,1]\rightarrow\R^+_0,\qquad \varphi(k):=\left(\gamma(\gamma-1)k^{1-\frac{2}{\gamma}}\right)^{-1}=\left(\gamma(\gamma-1)\right)^{-1}k^{\frac{2}{\gamma}-1}.
\end{align*}
We define a regularising transformation for $w$ by 
\begin{align}
 \Phi:[0,1]\rightarrow\R^+_0,\qquad\Phi(k):=\int_0^k\int_0^{s_2}\varphi(s_1)\,ds_1ds_2=\left(\gamma(\gamma-1)\frac{2}{\gamma}\left(\frac{2}{\gamma}+1\right)\right)^{-1}k^{\frac{2}{\gamma}+1}.\nonumber
 \end{align}
Then, we have that
\begin{align}
 \Phi(w(x))=\left(\gamma(\gamma-1)\frac{2}{\gamma}\left(\frac{2}{\gamma}+1\right)\right)^{-1}|x|^{2+\gamma}\text{ for all }x\in B_1.\nonumber
\end{align}
It is easy to see that $\Phi$ has the following properties:
\begin{enumerate}
  \item $\Phi(w)\in C^{2}(B_1)$;
 \item $\Phi$ is a strictly increasing function, so that it allows to reconstruct  back $w$ from $\Phi(w)$;
 \item $\Phi$ preserves the zero set of $w$  and allows to reconstruct it back;
  \item $\Phi$ can be used to regularise a whole class of functions $w$ which are  smooth everywhere but for their zero sets and satisfy \eqref{sing}. Thus, $\Phi$ smooths down such a function $w$ near its zero set. 
\end{enumerate}
\end{Example}\noindent
Let us now apply the idea from {\it Example \ref{simple}} to solutions $c$ of \eqref{c}. In this general case, however, we cannot hope for pointwise estimates  like \eqref{sing} to hold  uniformly in $\omega$. This is because the source term $f$ depends upon $y$, which is, for each $x$, a solution of an SDE. Hence, instead of using the spaces of functions which are differentiable in the classical sense, we  work in anisotropic Sobolev spaces. The, possibly, irregular behaviour of $c$ at the parabolic boundary of the cylinder $(0,T]\times O$ presents yet enough difficulty. Our regularising transformation should thus be able to smooth down $c$ not only near $\{c=0\}$, where the equation has a degeneracy, but also at the parabolic boundary  
\begin{align}
\Gamma:=\partial((0,T]\times O)\backslash (\{T\}\times O).\nonumber                                                                                                                                                                                                                                                                                                                                                                                                                                                                                                                                                                                                                                                                                                                                                             \end{align}
For this reason, we consider $c$ on the intersections of its level sets with a  decreasing family of subcylinders of $(0,T]\times O$:
\begin{align}
 Q_{d}:=\{(t,x)\in(0,T]\times O:\ \dist((t,x),\Gamma)> d\}\qquad\text{for all }d\in\left(0,\frac{1}{4}\diam(O)\right].\nonumber
\end{align}
We recall that due to  assumptions on $\beta$ and estimates \eqref{cbnd} and \eqref{fbnd} it holds (see, e.g., \cite[Theorem 2.I]{Ivanov86})  that
\begin{align}
&|c|_{C^{\alpha_0}({Q_d})}\leq \Cl{H1}\left(d^{-1}\right)\qquad\text{for all }d\in\left(0,\frac{1}{4}\diam(O)\right]\text{ for some }\alpha_0\in(0,1).\label{Hoelder}
\end{align}
Thus, the H\"older constant may explode as $d\rightarrow0$, that is, as $\Gamma$ is approached. 
We next divide both sides of \eqref{c} by $\beta'(c)$ and thus obtain an equation in  a non-divergence form:
\begin{align}
 &\partial_t c=(\beta'(c))^{-1}\Delta c+(\beta'(c))^{-1}f(c,y).\label{PM1}
\end{align}
 Due to assumption \eqref{betamm} and estimates \eqref{cbnd}, \eqref{fbnd}, and \eqref{Hoelder}, it holds 
that
\begin{align}
&\left\| (\beta'(c))^{-1}\right\|_{C^{\alpha_1}({Q_d})}\leq \Cl{H2}\left(d^{-1}\right)\text{ for some }\alpha_1\in(0,1),\label{H1}\\
&\left\|(\beta'(c))^{-1}f(c,y)\right\|_{L^{\infty}(Q_d)}\leq \C.\nonumber
\end{align}
In order to obtain \eqref{H1}, we used the well-known property of  superpositions of H\"older continuous functions: 
 \begin{align}
  &u_1\in C^{\gamma_1}({D}),\ u_2\in C^{\gamma_2}(u_1({D}))\text{ for some } \gamma_1,\gamma_2\in(0,1)\nonumber\\
  \Rightarrow &u_2\circ u_1\in C^{\gamma_1\gamma_2}({D})\text{ and }|u_2\circ u_1|_{C^{\gamma_1\gamma_2}({D})}\leq |u_2|_{C^{\gamma_2}(u_1({D}))}|u_1|_{C^{\gamma_1}({D})}^{\gamma_2}.\nonumber
 \end{align}
Let us consider for any $k\in\left(0,\left\|c\right\|_{L^{\infty}((0,T)\times O)}\right)$, $d\in\left(0,\frac{1}{4}\diam(O)\right]$ the sets $\{c>k\}\cap Q_{d}$  and $\left\{c>\frac{k}{2}\right\}\cap Q_{\frac{d}{2}}$. Using the crucial property \eqref{Hoelder}, we deduce the following: these two sets are relatively open w.r.t. $(0,T]\times O$ and their parabolic boundaries do not intersect. Moreover, the estimate on the H\"older norm allows to estimate the distance between the parabolic  boundaries from below by a positive number which depends only upon $d,k$, and, of course, the parameters of the problem. Now, equation \eqref{PM1} is non-degenerate in  $\left\{c>\frac{k}{2}\right\}$.   Therefore, we can apply standard results on local regularity for linear parabolic equations, see Theorems 9.1 and 10.1, and the remark on local estimates in Sobolev spaces at the end of \textsection 10 in \cite[Chapter IV]{LSU}.  Considering $\{c>k\}\cap Q_{d}$  as a subdomain of $\left\{c>\frac{k}{2}\right\}\cap Q_{\frac{d}{2}}$, these regularity results can be interpreted in the following way: for each $p\in (1,\infty)$  there exists a function
\begin{align}
 \varphi_p:\R^+_0\times\left[0,\frac{1}{4}\diam(O)\right] \rightarrow\R^+_0\nonumber
\end{align}
with the properties
\begin{enumerate}
\item\label{phi1}  $\varphi_p$ depends only upon $p$ and the parameters of the problem;
 \item\label{phi2}  $\varphi_p(0,\cdot)\equiv0$, $\varphi_p(\cdot,0)\equiv0$ in $\R^+_0$, $\varphi_p>0$ in $\R^+\times\left(0,\frac{1}{4}\diam(O)\right]$;
 \item\label{phi3}  $\varphi_p$ is uniformly bounded;
 \item\label{phi4}  $\varphi_p$ is increasing in each of the two variables;
 \item\label{phi5}  for each $k\in\left(0,\left\|c\right\|_{L^{\infty}((0,T)\times O)}\right)$ and $d\in\left(0,\frac{1}{4}\diam(O)\right]$ it holds that
 \begin{align}
 &\left\|c\right\|_{W^{(1,2),2p}(\{c>k\}\cap Q_d)}\leq \varphi_p^{-1}(k,d);\label{singc}
\end{align}
 \item\label{phi6} $\varphi_p\beta'\leq1$.
\end{enumerate}
\begin{Remark}
Local estimates from \cite{LSU} deal with cylindrical sets, which is generally not the case for a set  $\{c>k\}\cap Q_d$.  However, since the closure of such a set is a compact set, lies inside of a relatively w.r.t. $(0,T]\times O$ open set $\left\{c>\frac{k}{2}\right\}\cap Q_{\frac{d}{2}}$, and we have control upon the distance between the parabolic boundaries, we can cover $\{c>k\}\cap Q_d$ by a finite number of sufficiently small cylinders which lie completely in $\left\{c>\frac{k}{2}\right\}\cap Q_{\frac{d}{2}}$ and such that their number and the distance between their parabolic boundaries and the parabolic boundary of $\left\{c>\frac{k}{2}\right\}\cap Q_{\frac{d}{2}}$ is bounded from below by a positive number which depends only upon $d,k$, and the parameters of the problem. Hence, we can apply the results from \cite{LSU} to each of these cylinders and subsequently  sum together the resulting estimates in order to deduce \eqref{singc} with $\varphi$ satisfying conditions \ref{phi1}.-\ref{phi4}. from above.
  \end{Remark}
 \begin{Remark} Observe that if  a function $\tilde{\varphi}_p$ satisfies conditions \ref{phi1}.-\ref{phi5}. from above, we can clearly satisfy all six conditions by taking $$\varphi_p(k,d):=\min\{\tilde{\varphi}_p(k,d),(\beta'(k))^{-1}\}.$$
\end{Remark}
Estimate \eqref{singc} together with properties \ref{phi1}.-\ref{phi4}. convey that $c$ is well-behaved away from its zero set $\{c=0\}$ and the parabolic boundary of $(0,T]\times O$, may possibly have singularities on that parabolic boundary and/ or $\{c=0\}$, but, also, that we have some control on its behaviour near the singularities. 
Using $\varphi_p$, we are now able to produce our regularising transformation for $c$:
\begin{align}
 &\Phi_p:\R^+_0\times\left[0,\frac{1}{4}\diam(O)\right] \rightarrow\R^+_0,\nonumber\\
 &\Phi_p(k,d):=\int_0^k\int_0^N\int_0^{s_2}\int_0^{z_2}\frac{s_1}{2}\varphi_p^2\left(\frac{s_1}{2},z_1\right)\,dz_1\,ds_1\,dz_2\,ds_2.\nonumber
\end{align}
Due to properties \ref{phi2}. and \ref{phi4}. of $\varphi_p$, we have for each $k\in\left(0,\left\|c\right\|_{L^{\infty}((0,T)\times O)}\right)$ and $d\in\left(0,\frac{1}{4}\diam(O)\right]$  that 
\begin{align}
 0\leq \partial_{k^{\alpha_1}d^{\alpha_2}}^{\alpha_1+\alpha_2}\Phi_p(k,d)\leq \Cl{CPhi}\frac{k}{2}\varphi_p^2\left(\frac{k}{2},d\right) \text{ for }\alpha_1,\alpha_2\in\{0,1,2\},\ \alpha_1+\alpha_2\leq2. \label{Pribnd}
\end{align}
Next, we recall that domain $O$ has a smooth boundary. Consequently,  there exists a function 
\begin{align}
 \gamma:[0,T]\times \overline{O}\rightarrow\left[0,\frac{1}{4}\diam(O)\right]\nonumber
\end{align}
with the following properties:
\begin{enumerate}
 \item $\gamma(t,x)>0$ in $(0,T]\times O$, $\gamma(t,x)=0$ in $\Gamma$;
 \item $\gamma\in C^{2}([0,T]\times \overline{O})$;
 \item there exists a number $d_0\in \left(0,\frac{1}{4}\diam(O)\right]$, which depends only upon the domain $O$, such that
 \begin{align}
\gamma\leq d\text{ in } \overline{Q_d} \text{ for all }d\in[0,d_0].  \nonumber                                                                                              \end{align}
\end{enumerate}
Using the chain rule, \eqref{Pribnd} and the properties of $\varphi_p$ and $\gamma$, one readily checks that for all $p\in(1,\infty)$ it holds that
\begin{align}
 &\|\Phi_p(c,\gamma)\|_{W^{(1,2),p}((0,T)\times O)}\leq \C(p),\label{est4_pr}\\
 &\|\partial_t\Phi_p(c,\gamma)\beta'(c)\|_{L^p((0,T)\times O)}\leq \C(p).\label{est9}
\end{align}
Indeed, for each $i,j\in \{1,\dots,N\}$, $k\in\left(0,\left\|c\right\|_{L^{\infty}((0,T)\times O)}\right)$ and $d\in(0,d_0]$, it holds due to the H\"older inequality  that, for instance,
\begin{align}
\left\|\partial_{c^2}\Phi_p(c,\gamma)\partial_{x_i}c\partial_{x_j}c\right\|_{L^p(\{k<c\leq 2k\}\cap Q_d)}
\leq&\Cr{CPhi}\left\|\frac{c}{2}\varphi_p^2\left(\frac{c}{2},\gamma\right)\partial_{x_i}c\partial_{x_j}c \right\|_{L^p(\{k<c\leq 2k\}\cap Q_d)}\nonumber\\
\leq &\Cr{CPhi}k\varphi_p^2(k,d)\left\|\partial_{x_i}c\right\|_{L^{2p}(\{k<c\leq 2k\}\cap Q_d)}\left\|\partial_{x_i}c\right\|_{L^{2p}(\{k<c\leq 2k\}\cap Q_d)}\nonumber\\
\leq &\Cr{CPhi}k\varphi_p^2(k,d)\left\|c\right\|_{W^{(1,2),2p}(\{k<c\leq 2k\}\cap Q_d)}^2\nonumber\\
\leq &\Cr{CPhi}k.\label{indepk}
\end{align}
Since the constant $\Cr{CPhi}$ doesn't depend upon $k$, estimate \eqref{indepk} yields that
\begin{align}
&\left\|\partial_{c^2}\Phi_p(c,\gamma)\partial_{x_i}c\partial_{x_j}c\right\|_{L^p(\{c>0\}\cap Q_d)}\nonumber\\
=&\sum_{i=0}^{\infty}\left\|\partial_{c^2}\Phi_p(c,\gamma)\partial_{x_i}c\partial_{x_j}c\right\|_{L^p\left(\left\{2^{-(i+1)}\left\|c\right\|_{L^{\infty}((0,T)\times O)}<c\leq 2^{-i}\left\|c\right\|_{L^{\infty}((0,T)\times O)}\right\}\cap Q_d\right)}\nonumber\\
\leq&\Cr{CPhi}\sum_{i=0}^{\infty}2^{-(i+1)}\nonumber\\
=&\Cr{CPhi}.\label{est6_}
\end{align}
Finally, since $\Phi_p(c,\gamma)\equiv0$ on $\{c=0\}$ and $\Cr{CPhi}$ doesn't depend upon $d$, \eqref{est6_} implies that
\begin{align}
\left\|\partial_{c^2}\Phi_p(c,\gamma)\partial_{x_i}c\partial_{x_j}c\right\|_{L^p((0,T)\times O)}\leq\Cr{CPhi}.\nonumber
\end{align}
Applying the chain rule to $\Phi_p(c,\gamma)$ in order to compute the required partial derivatives and treating other resulting terms in a similar fashion leads to estimates \eqref{est4_pr} and \eqref{est9}.

\subsection{Estimates for a transformation of solutions of the RPME \texorpdfstring{\eqref{PM}}{}}\label{SecPsi}
Let $p\in(1,\infty)$ and let $\Phi_p$ be the smoothing transformation  from {\it Subsection \ref{secPMTr}}. 
We now introduce yet another transformation
\begin{align}
&\Psi_p: \R^+_0\times\left[0,\frac{1}{4}\diam(O)\right] \rightarrow\R^+_0,\nonumber\\
& \Psi_p(k,d):=\int_0^k\Phi_p(s,d)\beta'(s)\,ds.\nonumber
\end{align}
Due to assumption \eqref{betaderp} and the properties of $\Phi_p$, we have that 
\begin{align}
\Psi_p\in C^1\left(\R^+_0\times\left[0,\frac{1}{4}\diam(O)\right]\right),\qquad 0\leq\partial_k \Psi_p(k,d)\leq \Cl{CPsi0}(p,k). \label{DcPsi}
\end{align}
The continuity of $\partial_k \Psi_p$ in $\R^+\times\left[0,\frac{1}{4}\diam(O)\right]$ is a direct consequence of continuity of $\Phi_p$ and $\beta$. Moreover, it holds with \eqref{Pribnd} and the properties of $\varphi_p$ listed above that
\begin{align}
 0\leq\Phi_p(k,d)\beta'(k)\leq& \Cr{CPhi}\frac{k}{2}\varphi_p^2\left(\frac{k}{2},d\right)\beta'(k)\nonumber\\
 \leq&\Cl{CPsi1}(p)k\varphi_p\beta'(k)\nonumber\\
 \leq&\Cr{CPsi1}(p)k=:\Cr{CPsi0}(p,k).\label{est10}
\end{align}
Estimate \eqref{est10} yields in particular that $\partial_k \Psi_p$ is continuous in every point of the set $\{0\}\times \left[0,\frac{1}{4}\diam(O)\right]$, as required.

Using estimates  \eqref{cbnd}-\eqref{dtc} and  \eqref{DcPsi} and  the chain rule, we obtain the following group of estimates:
\begin{subequations}\label{Psipest}
\begin{empheq}[left={ \text{a.s. in }\Omega\ \empheqlbrace\,}]{align}
  &\left\|\Psi_p(c,\gamma)\right\|_{L^{\infty}((0,T);L^{\infty}(O))}\leq \Cl{Cphi}(p)\\
 &\left\|\nabla\Psi_p(c,\gamma)\right\|_{L^{\infty}((0,T);(L^2(O))^N)}\leq \Cr{Cphi}(p)\\
 &\left\|\partial_t \Psi_p(c,\gamma)\right\|_{L^2((0,T);L^2(O))}\leq \Cr{Cphi}(p),\label{cphi3}\\
 &\left\|D_r \Psi_p(c,\gamma)\right\|_{L^2((r,T)\times O)}\leq \Cr{Cphi}(p)\left\|D_r c\right\|_{L^2((r,T)\times O)},\label{cphi4}\\
 &\left\|D_{r+h_2}\Psi_p(c,\gamma)-D_r \Psi_p(c,\gamma)\right\|_{L^2((r,t)\times O)}
 \leq \Cr{Cphi}(p)\left\|D_{r+h_2}c-D_r c\right\|_{L^2((r,t)\times O)}. \label{cphi5}
\end{empheq}
\end{subequations}
Next, we use the chain and product rules in order to obtain the following representation for  $\partial_t D_r \Psi_p(c,\gamma)$:
\begin{align}
 \partial_t D_r \Psi_p(c,\gamma)=\Phi_p(c,\gamma)\partial_tD_r \beta(c)+\partial_t\Phi_p(c,\gamma)\beta'(c)D_r c.\label{psip}
 \end{align}
Let $p>\max\left\{\frac{N}{2},2\right\}$. Then, $W^{2,p}_0(O)$ is closed under pointwise multiplication and
\begin{align}
 \|uv\|_{W^{2,p}_0(O)}\leq\C(p)\|u\|_{W^{2,p}_0(O)}\|v\|_{W^{2,p}_0(O)},\nonumber
\end{align}
so that, due to the Sobolev inequality,
\begin{align}
 \|uv\|_{H^{2}_0(O)}\leq\Cl{C10}(p)\|u\|_{W^{2,p}_0(O)}\|v\|_{W^{2,p}_0(O)}.\label{algebra}
\end{align}
Using \eqref{algebra}, we obtain the following estimate for the first summand on the right-hand side of \eqref{psip}:
\begin{align}
 \left\|\Phi_p(c,\gamma)\partial_tD_r \beta(c)\right\|_{H^{-2}(O)}\leq \Cr{C10}(p)\left\|\Phi_p(c,\gamma)\right\|_{W^{2,p}(O)}\left\|\partial_tD_r \beta(c)\right\|_{H^{-2}(O)}.\label{1s}
\end{align}
For the second summand, the H\"older and Sobolev inequalities together with the properties of $\Phi_p$ yield that
\begin{align}
 \left\|\partial_t\Phi_p(c,\gamma)\beta'(c)D_r c\right\|_{H^{-2}(O)}\leq& \Cl{C11}(p)\left\|\partial_t\Phi_p(c,\gamma)\beta'(c)  D_r c\right\|_{L^{\frac{2p}{p+2}}(O)}\nonumber\\
 \leq &\Cr{C11}(p)\left\|\partial_t\Phi_p(c,\gamma)\beta'(c) \right\|_{L^p(O)}\left\|D_r c\right\|_{L^2(O)}.\label{2s}
\end{align}
Combining  \eqref{est4_pr}, \eqref{est9}, \eqref{psip}, \eqref{1s}, and \eqref{2s}, integrating over $(r,T)$ and using   the H\"older inequality, we obtain that
\begin{align}
\left\|\partial_t D_r \Psi_p(c,\gamma)\right\|_{L^{\frac{2p}{p+2}}(r,T;H^{-2}((r,T)\times O))}
\leq& \Cr{C10}(p)\left\|\Phi_p(c,\gamma)\right\|_{L^p(r,T;W^{2,p}(O))}\left\|\partial_tD_r \beta(c)\right\|_{L^2(r,T;H^{-2}(O))}\nonumber\\
&+\Cr{C11}(p)\left\|\partial_t\Phi_p(c,\gamma)\beta'(c) \right\|_{L^p((r,T)\times O)}\left\|D_r c\right\|_{L^2((r,T)\times O)}\nonumber\\
\leq&\Cl{C12}(p) \left(\left\|\partial_tD_r \beta(c)\right\|_{L^2(r,T;H^{-2}(O))}+\left\|D_r c\right\|_{L^2((r,T)\times O)}\right).\label{est7}
\end{align}
Since $p>2$, $C^{\frac{1}{2}-\frac{1}{p}}\subset W^{1,\frac{2p}{p+2}}$ holds. This, together with \eqref{cphi4} and \eqref{est7}, finally yields that 
\begin{align}
 \left\|D_r \Psi_p(c,\gamma)\right\|_{C^{\frac{1}{2}-\frac{1}{p}}([r,T];H^{-2}(O))}\leq&\Cl{C13}(p) \left(\left\|\partial_tD_r \beta(c)\right\|_{L^2(r,T;H^{-2}(O))}+\left\|D_r c\right\|_{L^2((r,T)\times O)}\right)\qquad \text{a.s. in }\Omega.\label{est8}
\end{align}
\section{A priori estimates for SDE \texorpdfstring{\eqref{SDE}}{}}\label{aprioriy}
\subsection{Basic estimate for \texorpdfstring{$y$}{y}}
We begin with an $L^q$-estimate  for $y$ as solution of the  stochastic integral equation \eqref{inteqy}. Thereby we choose $q>2$ as in the Sobolev embedding theorem, i.e., such that $H^1(O)\subset L^q(O)$. 
Using assumptions \eqref{aLip} and \eqref{bLip}, estimate \eqref{cbnd}, and a version of the Burkholder-Gundy-Davis inequality   \cite[Chapter 1 Theorem 7.1]{Ma}, we obtain  that
\begin{align}
 \|y(t)\|_{L^q(\Omega)}^q\leq &\Cl{C6}\left(|y_0|^q+\left\|\int_0^ta(y)(s)\,dW(s)\right \|_{L^q(\Omega)}^q+\left\|\int_0^tb(c,y)(s)\,ds\right\|^q_{L^q(\Omega)}\right)\nonumber\\
 \leq &\Cr{C6}|y_0|^q+t^{\frac{q-2}{2}}\Cr{C6}\int_0^t\|a(y)(s)\|_{L^q(\Omega)}^q\,ds+t^{q-1}\Cr{C6}\int_0^t\|b(c,y)(s)\|_{L^q(\Omega)}^q\,ds\nonumber\\
 \leq &\Cr{C6}|y_0|^q+\Cl{C32}\int_0^t\|y(s)\|_{L^q(\Omega)}^q\,ds+\Cr{C32}\int_0^tb^q(0,0)+\|c(s)\|_{L^q(\Omega)}^q+\|y(s)\|_{L^q(\Omega)}^q\,ds\nonumber\\
 \leq &\Cr{C6}|y_0|^q+\Cl{C27}+\Cr{C27}\int_0^t\|y(s)\|_{L^q(\Omega)}^q\,ds.\label{y210}
\end{align}
Integrating \eqref{y210} over $O$ and using assumption \eqref{y0}, we conclude that
\begin{align}
 \|y(t)\|_{L^q(\Omega;L^q(O))}^q
 \leq &\Cr{C6}\|y_0\|_{L^q(O)}^q+\C+\Cr{C27}\int_0^t\|y(s)\|_{L^q(\Omega;L^q(O))}^q\,ds\nonumber\\
 \leq &\Cl{C28}+\Cr{C27}\int_0^t\|y(s)\|_{L^q(\Omega;L^q(O))}^q\,ds.\label{y21}
\end{align}
Applying the Gronwall lemma to \eqref{y21}, we arrive at the estimate
\begin{align}
 \|y\|_{L^{\infty}((0,T);L^q(\Omega;L^q(O)))}\leq \Cl{C16}.\label{yLq}
 \end{align}
\subsection{Estimate for \texorpdfstring{$\nabla y$}{}}
Computing the spatial gradient on both sides of \eqref{inteqy}, we obtain that $\nabla y$ satisfies the stochastic integral equation
\begin{align}
 \nabla y(t)=\nabla y_0+\int_0^ta'(y)\nabla y(s)\,d W(s)+\int_0^t\partial_c b(c,y)\nabla c(s)+\partial_y b(c,y)\nabla y(s)\,ds.\label{nablaeqy}
\end{align}
Using assumptions \eqref{aLip} and \eqref{bLip}, estimates \eqref{cbnd} and \eqref{nablac}, and the It\^{o} isometry, we obtain   that
\begin{align}
 \|\nabla y(t)\|_{(L^2(\Omega))^N}^2
 \leq&4\|\nabla y_0\|_{(L^2(\Omega))^N}^2+4\int_0^t\|a'(y)\nabla y(s)\|_{(L^2(\Omega))^N}^2\,d s\nonumber\\
 &+4t\int_0^t\|\partial_c b(c,y)\nabla c(s)\|_{(L^2(\Omega))^N}^2+\|\partial_y b(c,y)\nabla y(s)\|_{(L^2(\Omega))^N}^2\,ds\nonumber\\
 \leq&4\|\nabla y_0\|_{(L^2(\Omega))^N}^2+\Cl{C30}\int_0^t\|\nabla c(s)\|_{(L^2(\Omega))^N}^2+\|\nabla y(s)\|_{(L^2(\Omega))^N}^2\,ds\nonumber\\
 \leq&4\|\nabla y_0\|_{(L^2(\Omega))^N}^2+\Cl{C31}+\Cr{C31}\int_0^t\|\nabla y(s)\|_{(L^2(\Omega))^N}^2\,ds.\label{nablay1}
\end{align}
Integrating \eqref{nablay1} over $O$ and using assumption \eqref{y0}, we conclude that
\begin{align}
 \|\nabla y(t)\|_{L^2(\Omega;(L^2(O))^N)}^2
\leq&4\|\nabla y_0\|_{L^2(O))^N}^2+\C+\Cr{C31}\int_0^t\|\nabla y(s)\|_{L^2(\Omega;(L^2(O))^N)}^2\,ds\nonumber\\
\leq&\C+\Cr{C31}\int_0^t\|\nabla y(s)\|_{L^2(\Omega;(L^2(O))^N)}^2\,ds.\label{naby21}
\end{align}
Applying the Gronwall lemma to \eqref{naby21}, we arrive at the estimate
\begin{align}
 \left\|\nabla y\right\|_{L^{\infty}((0,T);L^2(\Omega;(L^2(O))^N))}\leq \Cl{C41}.\label{nablay}
\end{align}
\subsection{Estimate for \texorpdfstring{$y(t+h_1)-y(t)$}{}}
The difference $y(t+h_1)-y(t)$ satisfies
\begin{align}
 y(t+h_1)-y(t)=\int_t^{t+h_1} a(y)(s)\,d W(s)+\int_t^{t+h_1} b(c,y)(s)\,ds.\nonumber
\end{align}
Using assumptions \eqref{aLip} and \eqref{bLip}, estimate \eqref{cbnd}, and a version of the Burkholder-Gundy-Davis inequality   \cite[Chapter 1 Theorem 7.1]{Ma}, we obtain that
\begin{align}
 &\left\|y(t+h_1)-y(t)\right\|_{L^q(\Omega)}^q\nonumber\\
 \leq& h_1^{\frac{q-2}{2}}\Cl{C100}\int_t^{t+h_1} \left\|a(y)(s)\right\|_{L^q(\Omega)}^q\,ds+h_1^{q-1}\Cr{C100} \int_t^{t+h_1} \left\|b(c,y)(s)\right\|_{L^q(\Omega)}^q\,d s\nonumber\\
 \leq& h_1^{\frac{q-2}{2}}\Cl{C32_}\int_t^{t+h_1}\|y(s)\|_{L^q(\Omega)}^q\,ds+h_1^{q-1}\Cr{C32_}\int_t^{t+h_1}b^q(0,0)+\|c(s)\|_{L^q(\Omega)}^q+\|y(s)\|_{L^q(\Omega)}^q\,ds\nonumber\\
 \leq&h_1^{q}\Cl{C33}+ h_1^{\frac{q-2}{2}}\Cr{C33}\int_t^{t+h_1}\|y(s)\|_{L^q(\Omega)}^q\,ds.\label{Dify1}
\end{align}
Integrating \eqref{Dify1} over $O$ and using  estimate \eqref{yLq}, we arrive at the estimate
\begin{align}
 \left\|y(t+h_1)-y(t)\right\|_{L^q(\Omega;L^q(O))}\leq h_1^{\frac{1}{2}}\Cl{C17}.\label{Difyhy}
\end{align}
\subsection{Estimate for the Malliavin derivative \texorpdfstring{$D_ry$}{Dry}}
Using the chain rule  and the rule of the differentiation of an It\^{o} SDE, we compute the $D_r$-derivative on both sides of equation \eqref{inteqy}. This leads to a stochastic integral equation for $D_r y$:
\begin{align}
 D_ry(t)=a(y)(r)+\int_r^ta'(y)D_ry(s)\,d W(s)+\int_r^t \partial_c b(c,y)D_rc(s)+ \partial_y b(c,y)D_ry(s)\,ds.\label{Dreqy_}
\end{align}
Using assumptions \eqref{aLip} and \eqref{bLip}, estimate \eqref{cbnd}, and the It\^{o} isometry, we obtain  that
\begin{align}
 \|D_r y(t)\|_{L^2(\Omega)}^2
 \leq&4\|a(y)(r)\|_{L^2(\Omega)}^2+4\int_r^t\|a'(y)D_r y(s)\|_{L^2(\Omega)}^2\,d s\nonumber\\&+4(t-r)\int_r^t\|\partial_c b(c,y)D_r c(s)\|_{L^2(\Omega)}^2+\|\partial_y b(c,y)D_r y(s)\|_{L^2(\Omega)}^2\,ds\nonumber\\
 \leq&\Cl{C34}\|y(r)\|_{L^2(\Omega)}^2+\Cr{C34}\int_r^t\|D_r y(s)\|_{L^2(\Omega)}^2\,ds+\Cr{C34}\int_r^t\|D_r c(s)\|_{L^2(\Omega)}^2\,ds.\label{Dry1_}
\end{align}
Integrating \eqref{Dry1_} over $O$ and using estimate \eqref{yLq}, we conclude that
\begin{align}
 \|D_r y(t)\|_{L^2(\Omega;L^2(O))}^2
 \leq&\Cr{C34}\|y(r)\|_{L^2(\Omega;L^2(O))}^2+\Cr{C34}\int_r^t\|D_r y(s)\|_{L^2(\Omega;L^2(O))}^2\,ds\nonumber\\
 &+\Cr{C34}\int_r^t\|D_r c(s)\|_{L^2(\Omega;L^2(O))}^2\,ds\nonumber\\
 \leq&\Cl{C35}+\Cr{C34}\int_r^t\|D_r y(s)\|_{L^2(\Omega;L^2(O))}^2\,ds+\Cr{C34}\int_r^t\|D_r c(s)\|_{L^2(\Omega;L^2(O))}^2\,ds.\label{Dry2_}
\end{align}
Applying the Gronwall lemma to \eqref{Dry2_}, we arrive at the estimate
\begin{align}
\left\|D_ry(t)\right\|_{L^2(\Omega;L^2(O))}^2\leq& \Cl{C5}+\Cr{C5}\int_r^t\left\|D_r c(s)\right\|_{L^2(\Omega;L^2(O))}^2\,ds.\label{Dry1}
\end{align}
 \subsection{Estimate for  \texorpdfstring{$D_r y(t+h_1)-D_r y(t)$}{}}
 Due to \eqref{Dreqy_}, the difference $D_r y(t+h_1)-D_r y(t)$ satisfies 
\begin{align}
 D_r y(t+h_1)-D_r y(t)=&\int_t^{t+h_1}a'(y(s))D_ry(s)\,d W(s)+\int_t^{t+h_1} \partial_c b(c,y)D_rc(s)+ \partial_y b(c,y)D_ry(s)\,ds.\nonumber
\end{align}
Using assumptions \eqref{aLip} and \eqref{bLip}, estimate \eqref{cbnd}, and the It\^{o} isometry, we obtain  that
\begin{align}
 &\left\|D_r y(t+h_1)-D_r y(t)\right\|_{L^2(\Omega)}^2\nonumber\\
 \leq&3\int_t^{t+h_1}\left\|a'(y(s))D_ry(s)\right\|_{L^2(\Omega)}^2\,d s+3h_1\int_t^{t+h_1} \left\|\partial_c b(c,y)D_rc(s)\right\|_{L^2(\Omega)}^2+ \left\|\partial_y b(c,y)D_ry(s)\right\|_{L^2(\Omega)}^2\,ds\nonumber\\
 \leq&\Cl{C38_}\int_t^{t+h_1}\left\|D_ry(s)\right\|_{L^2(\Omega)}^2\,d s+h_1\Cr{C38_}\int_t^{t+h_1}\left\|D_rc(s)\right\|_{L^2(\Omega)}^2\,ds.\label{Dreqy}
\end{align}
Integrating \eqref{Dry1_} over $O$ and using estimate \eqref{Dry1}, we conclude that
\begin{align}
 \left\|D_r y(t+h_1)-D_r y(t)\right\|_{L^2(\Omega;L^2(O))}^2
 \leq&h_1\Cr{C15}+h_1\Cl{C15}\int_r^{t+h_1}\left\|D_r c(s)\right\|_{L^2(\Omega;L^2(O))}^2\,ds.\label{Dreqy__}
\end{align}
 \subsection{Estimate for  \texorpdfstring{$D_{r+h_2}y-D_ry$}{}}
Due to \eqref{Dreqy_}, the difference $D_{r+h_2}y-D_ry$ satisfies (recall that $D_{t_1} F(t_2)=0$ for $t_2<t_1$ for $F$ adapted)
 \begin{align}
 D_{r+h_2}y(t)-D_r y(t)=&  a(y)(r+h_2)-a(y)(r)+\int_r^ta'(y) (D_{r+h_2}y-D_r y)(s)\,d W(s)\nonumber\\
 &+\int_r^t \partial_c b(c,y) (D_{r+h_2}c-D_r c)(s)+\partial_y b(c,y) (D_{r+h_2}y-D_r y)(s)\,ds
\end{align}
Using assumptions \eqref{aLip} and \eqref{bLip}, estimate \eqref{cbnd}, and the It\^{o} isometry, we obtain  that
\begin{align}
 &\left\|(D_{r+h_2}y-D_r y)(t)\right\|_{L^2(\Omega)}^2\nonumber\\
 \leq & 4\|a(y)(r+h_2)-a(y)(r)\|^2+4\int_r^t\|a'(y) (D_{r+h_2}y-D_r y)(s)\|_{L^2(\Omega)}^2\,d s\nonumber\\
 &+4(t-r)\int_r^t \|\partial_c b(c,y) (D_{r+h_2}c-D_r c)(s)\|_{L^2(\Omega)}^2+\|\partial_y b(c,y) (D_{r+h_2}y-D_r y)(s)\|_{L^2(\Omega)}^2\,ds\nonumber\\
 \leq & \Cl{C39}\|y(r+h_2)-y(r)\|_{L^2(\Omega)}^2+\Cr{C39}\int_r^t \|(D_{r+h_2}c-D_r c)(s)\|_{L^2(\Omega)}^2+\|(D_{r+h_2}y-D_r y)(s)\|_{L^2(\Omega)}^2\,ds.\label{Drh21_}
\end{align}
Integrating \eqref{Drh21_} over $O$ and using estimate \eqref{Dry1}, we conclude that
\begin{align}
 &\left\|(D_{r+h_2}y-D_r y)(t)\right\|_{L^2(\Omega;L^2(O))}^2\nonumber\\
\leq & \Cr{C39}\|y(r+h_2)-y(r)\|_{L^2(\Omega;L^2(O))}^2\nonumber\\
&+\Cr{C39}\int_r^t \|(D_{r+h_2}c-D_r c)(s)\|_{L^2(\Omega;L^2(O))}^2+\|(D_{r+h_2}y-D_r y)(s)\|_{L^2(\Omega;L^2(O))}^2\,ds\nonumber\\
\leq & h_2\Cl{C40}+\Cr{C39}\int_r^t \|(D_{r+h_2}c-D_r c)(s)\|_{L^2(\Omega;L^2(O))}^2\,ds+\int_r^t \|(D_{r+h_2}y-D_r y)(s)\|_{L^2(\Omega;L^2(O))}^2\,ds.\label{Drh21}
\end{align}
Applying the Gronwall lemma to \eqref{Drh21}, we arrive at the estimate
\begin{align}
\left\|(D_{r+h_2}y-D_r y)(t)\right\|_{L^2(\Omega;L^2(O))}^2\leq h_2\Cl{C8}
+\Cr{C8}\int_r^t\left\|(D_{r+h_2}c-D_r c)(s)\right\|_{L^2(\Omega;L^2(O))}^2\,ds.\label{DrDry2}
\end{align}
\section{Proof of the compactness {\it Theorem \ref{CompTh}}}\label{SecComp}
In this section, we finally prove our main result, {\it Theorem \ref{CompTh}} on compactness. We begin with collecting together estimates \eqref{cbnd}-\eqref{dtc}, \eqref{estDrc01}, \eqref{estDrc4}, \eqref{estDrbeta}, \eqref{nablay}, \eqref{yLq}, \eqref{Difyhy}, \eqref{Dry1}, \eqref{DrDry2}, \eqref{Dreqy__} for $c$ and $y$ which we obtained in {\it Sections \ref{aprioric}} and {\it \ref{aprioriy}}: 
\begin{subequations}\label{estc}
\begin{empheq}[left={\text{a.s. in } \Omega\ \empheqlbrace\,}]{align}
&\left\|c\right\|_{L^{\infty}((0,T);L^{\infty}(O))}\leq \Cr{Cap},\\
 &\left\|\nabla c\right\|_{L^{\infty}((0,T);(L^2(O))^N)}\leq \Cr{Cap},\\
 &\left\|\partial_t c\right\|_{L^2((0,T);L^2(O))}\leq \Cr{Cap},\\
&\int_r^t\left\|D_rc(s)\right\|_{L^2(O)}^2\,ds\leq 
\Cr{CDr1}\int_r^t\left\|D_ry(s)\right\|_{L^2(O)}^2\,ds,\\
& \int_r^t\left\|(D_{r+h_2}c-D_r c)(s)\right\|_{L^2(O)}^2\,ds
 \leq \Cr{CDr2}\int_r^t\left\|(D_{r+h_2}y-D_r y)(s)\right\|_{L^2(O)}^2\,ds,\\
 &\int_r^t\left\|\partial_t D_r\beta(c)(s)\right\|_{H^{-2}(O)}^2\,ds\leq \Cr{CdtDr1}\int_r^t\left\|D_ry(s)\right\|_{L^2(O)}^2\,ds
\end{empheq}
and
\end{subequations}
\begin{subequations}\label{esty}
\begin{align}
&\|y\|_{C^{\frac{1}{2}}([0,T];L^q(\Omega;L^q(O)))}\leq\Cr{C17},\label{yth1}\\
&\|\nabla y\|_{L^{\infty}((0,T);L^2(\Omega;(L^2(O))^N))}\leq \Cr{C41},\\
&\left\|D_ry(t)\right\|_{L^2(\Omega;L^2(O))}^2\leq  \Cr{C5}+\Cr{C5}\int_r^t\left\|D_r c(s)\right\|_{L^2(\Omega;L^2(O))}^2\,ds,\\
&\left\|(D_{r+h_2}y-D_r y)(t)\right\|_{L^2(\Omega;L^2(O))}^2\leq h_2\Cr{C8}
+\Cr{C8}\int_r^t\left\|(D_{r+h_2}c-D_r c)(s)\right\|_{L^2(\Omega;L^2(O))}^2\,ds,\\
&\left\|D_r y(t+h_1)-D_r y(t)\right\|_{L^2(\Omega;L^2(O))}^2
 \leq h_1\Cr{C15}+h_1\Cr{C15}\int_r^{t+h_1}\left\|D_r c(s)\right\|_{L^2(\Omega;L^2(O))}^2\,ds.\label{Drth1y}
\end{align}
\end{subequations}
Combining  \eqref{estc} and \eqref{esty} and using the Gronwall lemma where necessary, we arrive at the following set of estimates:
\begin{subequations}\label{estc1}
\begin{align}
&\left\|c\right\|_{L^{\infty}(\Omega;L^{\infty}((0,T);L^{\infty}(O)))}\leq \Cr{Cap},\\
 &\left\|\nabla c\right\|_{L^{\infty}(\Omega;L^{\infty}((0,T);(L^2(O))^N))}\leq \Cr{Cap},\\
 &\left\|\partial_t c\right\|_{L^{\infty}(\Omega;L^2((0,T);L^2(O)))}\leq \Cr{Cap},\\
&\left\|D_rc\right\|_{L^2((r,T)\times O\times\Omega)}\leq 
\Cl{tot1},\\
&\left\|h_2^{-\frac{1}{2}} \left(D_{r+h_2}c-D_r c\right)\right\|_{L^2((r,T-h_2)\times O\times \Omega)}
 \leq \Cr{tot1},\\
 &\left\|\partial_t D_r\beta(c)\right\|_{L^2((r,T)\times\Omega;H^{-2}(O))}\leq \Cr{tot1}\label{estDrc3__}
\end{align}
and
\end{subequations}
\begin{subequations}\label{esty1}
\begin{align}
&\|y\|_{C^{\frac{1}{2}}([0,T];L^q(\Omega;L^q(O)))}\leq\Cr{C17},\\
&\|\nabla y\|_{L^{\infty}((0,T);L^2(\Omega;(L^2(O))^N))}\leq \Cr{C41},\\
&\left\|D_ry\right\|_{L^{\infty}(r,T;L^2(\Omega;L^2(O)))}\leq\Cr{tot1},\\
&\left\|h_2^{-\frac{1}{2}}(D_{r+h_2}y-D_r y)\right\|_{L^{\infty}(r,T-h_2;L^2(\Omega;L^2(O)))}\leq\Cr{tot1}\\
&\left\|h_1^{-\frac{1}{2}}(D_r y(\cdot+h_1)-D_r y)\right\|_{L^{\infty}(r,T-h_1;L^2(\Omega;L^2(O)))}
 \leq\Cr{tot1}.
\end{align}
\end{subequations}
Estimates \eqref{esty1} allow us to apply \cite[Theorem 2]{BallySaussereau} directly. It  yields that
\begin{align}
&{\mathcal{Y}}(T,R)\text{ is precompact in } L^2((0,T)\times O\times\Omega).\nonumber
 \end{align}
 \begin{Remark}
 Observe that our estimates, particularly those  involving $D_r$, are in fact stronger than those required by that theorem. Indeed, for instance, assumptions (2)-(4) from \cite[Theorem 2]{BallySaussereau} deal with the regularised versions of the functions in the family.
 \end{Remark}
 Let us know prove the precompactness of ${\mathcal{C}}(T,R)$. First, we note that  \eqref{estDrc3__} is  an estimate for a second derivative of  $\beta(c)$, not for $c$.
This precludes the direct application of  \cite[Theorem 2]{BallySaussereau}. To overcome this problem, we consider instead function    $\Psi_p(c,\gamma)$. Combining the estimates \eqref{Psipest} and \eqref{est8}, which we derived in {\it Subsection \ref{SecPsi}}, with \eqref{estc1}, we obtain that
\begin{subequations}\label{estpsi}
 \begin{align}
 &\left\|\Psi_p(c,\gamma)\right\|_{L^{\infty}(\Omega;L^{\infty}((0,T);L^{\infty}(O)))}\leq \Cr{Cphi}(p),
 \\
  &\left\|\nabla \Psi_p(c,\gamma)\right\|_{L^{\infty}(\Omega;L^{\infty}((0,T);(L^2(O))^N))}\leq \Cr{Cphi}(p),
 \\
 &\left\|\partial_t \Psi_p(c,\gamma)\right\|_{L^{\infty}(\Omega;L^2((0,T);L^2(O)))}\leq \Cr{Cphi}(p),
 \\
 &\left\|D_r \Psi_p(c,\gamma)\right\|_{L^2((r,T)\times O\times\Omega)}\leq \Cl{C13_}(p),
 \\
 &\left\|h_2^{-\frac{1}{2}}\left(D_{r+h_2}\Psi_p(c,\gamma)-D_r \Psi_p(c,\gamma)\right)\right\|_{L^2((r,T-h_2)\times O\times\Omega)}
 \leq \Cr{C13_}(p), 
 \\
 &\left\|D_r \Psi_p(c,\gamma)\right\|_{L^2(\Omega;C^{\frac{1}{2}-\frac{1}{p}}([r,T];H^{-2}(O)))}\leq\Cr{C13_}(p).
\end{align}
\end{subequations}
With estimates \eqref{estpsi} at hand we can  know apply \cite[Theorem 2]{BallySaussereau} yielding that
\begin{align}
  &\Psi_p({\mathcal{C}}(T,R),\gamma)\text{ is precompact in } L^2((0,T)\times O\times\Omega).\label{compphi}
 \end{align}
 We observe that function $\Psi_p(\cdot,\gamma(t,x))$ has, for each fixed pair $(t,x)\in (0,T]\times O$, the following properties: it is defined on an interval, continuous, and strictly increasing. As to the latter, it is follows from the definition of $\Psi_p$ and the fact that $\beta'(c),\Phi_p(c,d)>0$ for $c,d>0$ and $\gamma(t,x)>0$ for $(t,x)\in (0,T]\times O$. 
 Therefore, $\Psi_p(\cdot,\gamma(t,x))$ is invertible, and its inverse has these three properties, too. Consequently, we have the following implication:
 \begin{align}
  &\left\{\Psi_p(c_n(t,x,\omega),\gamma(t,x))\right\}_{n\in\N}\text{ is convergent in }\R\Rightarrow
  \{c_n(t,x,\omega)\}_{n\in\N}\text{ is convergent in }\R.\label{nim}
 \end{align}
Combining \eqref{cbnd} and \eqref{nim} and using the dominated convergence theorem, we obtain that
 \begin{align}
  &\left\{\Psi_p(c_n,\gamma)\right\}_{n\in\N}\text{ is a.e. convergent in }(0,T)\times O\times\Omega
  \Rightarrow \{c_n\}_{n\in\N}\text{ is convergent in }L^2((0,T)\times O\times\Omega).\label{nim2}
 \end{align}
Together, \eqref{compphi} and  \eqref{nim2} finally yield  that
\begin{align}
 &{\mathcal{C}}(T,R)\text{ is precompact in } L^2((0,T)\times O\times\Omega).\nonumber
\end{align}
The proof of {\it Theorem \ref{CompTh}} is thus complete.

\section{Spatial semi-discretization for a nondegenerate case}\label{SecSemi}
In this section we set up and study a spatial finite-difference scheme for system \eqref{PM}-\eqref{SDE} under the additional 
 assumption
\begin{align}
\beta\in C^2(\R_0^+)\label{nondeg}
\end{align}
which corresponds to a nondegenerate case.
Our goal here is twofold. On one hand, we apply the semi-discretisation method (see, e.g. \cite[Chapter 4 \S 1]{Lions} and references therein, particularly \cite{raviart1970}) in order to obtain the existence of solutions. At the same time, we illustrate thereby how one can use compactness in order to rigorously prove the convergence of a numerical scheme for nonlinear systems with stochasticity. 


\subsection{Spatial discretisation and interpolation}
In this sequel we recall some concepts and ideas of the deterministic semi-discretisation method. We refrain from the  proofs of the properties listed here since they either exactly repeat or are slight modifications of results addressed in the literature. The interested reader is referred  to  \cite[Chapter 4 \S 1]{Lions}, where the method is described.  
The discretization is preformed only in $O$. To avoid some purely deterministic technical  difficulties which have to do with discretizing close to and on the boundary of $O$,   we restrict our exposition to the case when the spatial domain $O$ is a unit cube:
$$O=(0,1)^N.$$ However, we especially emphasise that this simplification is by no means essential: the present approach can be used for more general domains. 

We begin with some more notation. As usual, we denote by $e_k$ the $k$-th standard basis vector in $\R^N$.  Let  $M\in\N\backslash\{1\}$. For $h:=\frac{1}{M+1}$ we define  the discrete sets
\begin{align*}
 &O^h:=\left\{h,2h,\dots,1-h\right\}^N,\qquad \overline{O}^h:=\left\{0,h,\dots,1\right\}^N,\qquad \partial^hO^h:=\overline{O}^h\backslash O^h
\end{align*}
and the space-continuous set
\begin{align}
  &U^h:=(2h,1-2h)^N.\nonumber
\end{align}
Clearly, it holds that
\begin{align}
&|O\backslash U^h|=1-(1-4h)^N\underset{h\rightarrow0}{\rightarrow}0.\nonumber
\end{align}
Next, we introduce standard  finite difference operators for a function $u$ at a point $m$:
\begin{align}
 &\partial_k^hu(m):=\frac{1}{h}(u(m+he_k)-u(m))\quad\text{ for }
 k\in \{1,\dots,N\},\label{dkh}\\
 &\nabla^h:=\left(\partial_1^h,\dots,\partial_N^h\right)^T,\label{nablah}\\
 &\Delta^hu(m):=\frac{1}{h^2}\sum_{k=1}^{N}\left(u(m+he_k)-2u(m)+u(m-he_k)\right),
 \\
 &\partial_{\nu}^hu(m):=
   \frac{1}{h}(u(m)-u(\overline{m}))\quad\text{ for }m\in \partial^hO^h
\end{align}
where
\begin{align}
 \overline{m}_k:=\begin{cases}
   h&\text{ if }m_k=0,\\
   m_k&\text{ if }m_k\in(0,1),\\
   1-h&\text{ if }m_k=1.                                                                       \end{cases}\quad\text{ for } k\in\{1,\dots,N\}.
   \label{nuh}
\end{align}
We recall the discrete version of the  Green's first identity (i.e., summation by parts formula)
 \begin{align}
 \sum_{i=1}^M(a_{i+1}-2a_i+a_{i-1})b_i=-\sum_{i=0}^M(a_{i+1}-a_i)(b_{i+1}-b_i)+(a_{M+1}-a_M)b_{M+1}-(a_1-a_0)b_0,\label{GrI}
\end{align}
as well as the following relation 
\begin{align}
 \partial_k^h(g(u))(m)=\partial_k^hu(m)\int_0^1g'\left(u(m)+\tau h \partial_k^hu(m)\right)\,d\tau\label{chain}
\end{align}
which serves as  a sort of chain rule for the discrete case.
Further, we make use of some discrete  analogs of several  (semi-)norms and a scalar product which mimic the corresponding notions in connection with the  Lebesgue and Sobolev spaces: for $u,v:\overline{O}^h\rightarrow\R$ let
\begin{align}
 &(u,v)_{L^2(D^h;h)}:=h^{N}\sum_{m\in D^h}uv(m)\quad\text{ for }D^h\subset \overline{O}^h,\nonumber\\
 &\|u\|_{L^p(D^h;h)}:=\begin{cases}
\left(h^N\sum_{m\in D^h}|u(m)|^p\right)^{\frac{1}{p}}&\text{for } p\in[1,\infty),\\
\max_{m\in D^h}|u(m)|&\text{for }  p=\infty,
 \end{cases}\quad\text{ for }D^h\subset \overline{O}^h,\nonumber\\
 &|u|_{H^1(\overline{O}^h;h)}:=\left(h^{N-2}\sum_{m,m+he_k\in \overline{O}^h}|u(m+he_k)-u(m)|^2\right)^{\frac{1}{2}},\nonumber\\
 &\|u\|_{H^1(\overline{O}^h;h)}:=\left(|u|_{H^1(\overline{O}^h;h)}^2+\|u\|_{L^2(\overline{O}^h;h)}^2\right)^{\frac{1}{2}}\nonumber\\
 &\|u\|_{H^{-2}(\overline{O}^h;h)}:=\sup_{v\in H^2_0(\overline{O}^h;h)\backslash\{0\}}\frac{(u,v)_{L^2(O^h;h)}}{\|\Delta^h v\|_{L^2(O^h;h)}},\nonumber
\end{align}
where
\begin{align}
H^2_0(\overline{O}^h;h):=\left\{v:\overline{O}^h\rightarrow\R:\ v=\partial_{\nu}^hv=0\text{ on }\partial^hO^h \right\}.\nonumber
\end{align}
For a function $u:\overline{O}\rightarrow\R$ we define its projection to the space of discrete functions via
\begin{align}
P^h u(m):=\avint_{\left(m+\left(-\frac{h}{2},\frac{h}{2}\right)^N\right)\cap O} u(z)\,dz\quad\text{ for }m\in\overline{O}^h.\nonumber
\end{align}
In order to interpolate  discrete functions $u:\overline{O}^h\rightarrow \R$ 
 we use two types of splines: the  piecewise constant  
\begin{align}
 &\Pi^hu(x):=\begin{cases}
 u(m)&\text{ for }x\in \left(m+\left(-\frac{h}{2},\frac{h}{2}\right)^N\right)\cap \overline{O},\ m\in \overline{O}_h,\\
 0&\text{ otherwise}
\end{cases}
\nonumber
\end{align} 
 and the piecewise polyaffine
\begin{align}
 &\Lambda^hu(x):=\underset{i\in \left(m+\{0,h\}^N\right)\cap \overline{O}_h}{\sum}u(i)L_{\frac{1}{h}\left(i-m\right)}\left(\frac{1}{h}(x-m)\right)\quad\text{ for }x\in \left(m+[0,h]^N\right)\cap \overline{O},\ m\in \overline{O}_h,\nonumber
\end{align}
where
\begin{align}
&\Delta^h(z):=\prod\limits_{j=1}^N\beta_{l_j}(z_j)\quad \text{ for }z\in[0,1]^N,\ l\in \{0,1\}^N,\nonumber\\
  &\beta_0(r):=1-r,\quad \beta_1(r):=r\quad \text{ for }r\in[0,1].\nonumber
\end{align}
The constructed splines have the following useful properties: for all $u:\overline{O}^h\rightarrow \R$ it holds that  
\begin{align}
&\Pi^h u\in L^{\infty}(O),\nonumber\\
&(\Pi^h u,\Pi^h v)_{L^2(O)}=(u,v)_{L^2(\overline{O}^h;h)}\quad\text{ for }v\in H^2_0(\overline{O}^h;h)\label{lift},\\
&\varphi(\Pi^h u)=\Pi^h \varphi(u) \quad\text{ a.e. in }O\text{ for any } \varphi:\R\rightarrow\R, \label{PiProp}
\end{align}
\begin{align}
&\Lambda^hu\in W^{1,\infty}(O),\nonumber\\                                                     
 &\Lambda^hu(m)=u(m)\quad\text{ for } m\in \overline{O}^h,\nonumber
\end{align} 
and 
\begin{alignat}{4}
 u=0\quad\text{ in }\partial^hO^h&\quad\Rightarrow\quad&&\Lambda^hu=0 \quad&\text{ in }\partial O,\label{Lambda0}\\
 \partial_{\nu}^hu=0\quad\text{ in }\partial^hO^h&\quad\Rightarrow\quad&&\partial_{\nu}\Lambda^hu=0 \quad&\text{ a.e. in }\partial O.\label{nLambda0}
\end{alignat}
Moreover, both interpolation operators preserve positivity:
\begin{align}
 u\geq0\quad\text{ in }\overline{O}^h\quad\Rightarrow\quad \Lambda^hu,\Pi^hu\geq0\quad\text{ in }\overline{O}.\label{LPpos}
\end{align}
Finally, the projection and interpolations enjoy the following estimates:
\begin{alignat}{3}
 &\left\|P^hu\right\|_{L^p(\overline{O}^h;h)}\leq \Cl{CProj}\|u\|_{L^p(O)}&&\quad\text{ for }u\in L^p(O),\ p\in[1,\infty], \label{EstPr1}\\
 &\left|P^h u\right|_{H^1(\overline{O}^h;h)}\leq \Cr{CProj}\|\nabla u\|_{L^2(O)}&&\quad\text{ for }u\in H^1(O), \label{EstPr2}\\
&\left\|\Pi^hP^h u-u\right\|_{L^2(O)}\leq h\Cr{CProj}\|\nabla u\|_{L^2(O)}&&\quad\text{ for }u\in H^1(O),\label{estPiP}\\
&\left\|\Lambda^hu\right\|_{L^p(O)}\leq \Cr{CProj}\|u\|_{L^p(\overline{O}^h;h)}&&\quad\text{ for }u:\overline{O}^h\rightarrow \R,\  p\in[1,\infty],\label{TransL1}\\
 &\left\|\nabla\Lambda^hu\right\|_{L^2(O)}\leq \Cr{CProj}|u|_{H^1(\overline{O}^h;h)}&&\quad\text{ for }u:\overline{O}^h\rightarrow \R,\\
 &\left\|\Lambda^hu\right\|_{H^{-2}(U^h)}\leq \Cr{CProj}\left\|u\right\|_{H^{-2}(\overline{O}^h;h)}&&\quad\text{ for }u:\overline{O}^h\rightarrow \R,\label{TransL3}\\
 &\left\|\Pi^hu-\Lambda^hu\right\|_{L^2(O)}\leq h\Cr{CProj}\left|\nabla^h u\right|_{H^1(\overline{O}^h;h)}&&\quad\text{ for }u:\overline{O}^h\rightarrow \R,\label{compare}\\
 &\left\|\Pi^{h}\Delta^{h}u-u\right\|_{L^{\infty}(O)}\leq h\Cr{CProj}&&\quad\text{ for }u\in W^{3,\infty}(O).\label{PiDelh}
\end{alignat}

\subsection{Approximation via semi-discretization}
We start with the following  semi-discretization of the RPDE-SDE system \eqref{PM}-\eqref{SDE}:
\begin{subequations}\label{PMdiscr}
 \begin{empheq}[left={\text{in } \Omega\ \empheqlbrace\,}]{align}
&\partial_t \beta\left(c^h\right)=\Delta^h c^h+f\left(c^h,y^h\right)&&\text{ in }(0,T]\times O^h,\label{cdiscr}\\
 &c^h=0\qquad (\partial_{\nu}^hc^h=0)&&\text{ in }(0,T]\times \partial^hO^h,\label{bcdiscr}\\
 &c^h=c^h_0&&\text{ in }\{0\}\times O^h,\label{inic0discr}
\end{empheq}
\end{subequations}
\begin{subequations}\label{SDEdiscr}
 \begin{empheq}[left={\text{in } \overline{O}^h\ \empheqlbrace\,}]{align}
 &dy^h=a\left(y^h\right)\,d W+b\left(c^h,y^h\right)\,dt &&\text{ in } (0,T]\times\Omega,\label{ydiscr}\\
 &y^h=y^h_0&&\text{ in } \{0\}\times\Omega\label{iniy0discr}
\end{empheq}
\end{subequations}
where the spatially discretized  initial data $\left(c^h_0,y^h_0\right):\overline{O}^h\rightarrow\R_0^+$ is defined via
\begin{align}
\left(c^h_0,y^h_0\right):=\left(P^hc_0,P^hy_0\right)
\quad\text{ in }\overline{O}^h.\nonumber
\end{align}
Observe that equation \eqref{cdiscr} can be rewritten in the conventional form of a RODE:
\begin{align}
 \partial_t c^h =\left(\beta'\left(c^h\right)\right)^{-1}\Delta^h c^h+\left(\beta'\left(c^h\right)\right)^{-1}f\left(c^h,y^h\right).\label{conven}
\end{align}
 Using the boundary conditions \eqref{bcdiscr} one can eliminate  $\left(c^h(\cdot,m,\cdot)\right)_{m\in \partial^hO^h}$ from  \eqref{cdiscr}. Then, system \eqref{PMdiscr}-\eqref{SDEdiscr} can be considered as a RODE-SDE system  with respect to $$\left(\left(c^h(\cdot,m,\cdot)\right)_{m\in O^h},\left(y^h(\cdot,m,\cdot)\right)_{m\in \overline{O}^h}\right):[0,T]\times\Omega\rightarrow (\R_0^+)^{M^N}\times (\R_0^+)^{(M+2)^N}.$$ 
 A solution  can then be understood in the usual strong SDE-sense. Observe that not only  the coefficient functions $f,a,b$, but also $(\beta')^{-1}$ and $\Delta^h:\overline{O}^h\rightarrow O^h$ in equations \eqref{conven} and \eqref{ydiscr} are continuously differentiable in $\R^+_0$. 
 For $\beta$ this holds since $(\beta')^{-1},\beta''\in C(\R^+_0)$, the latter due to our additional assumption \eqref{nondeg}, while $\Delta^h$ is simply linear. Therefore, the RODE-SDE is uniquely solvable and its solutions possess square integrable  Malliavin derivatives due to well-known results \cite[Corollary 2.2.1, Theorem 2.2.1]{nualart2006malliavin}.  The nonnegativity of solution components is a consequence of a general result on the  invariance for SDE systems with smooth coefficients \cite{Milian} and assumptions \eqref{pospar} and \eqref{inipos}.

Multiplying \eqref{cdiscr} by $\xi v^h$ for  arbitrary  $\xi\in C_0^{\infty}[0,T)$ and $v^h\in H^2_0(\overline{O}^h;h)$, summing over $O^h$ thereby using \eqref{GrI} twice, and  integrating over $(0,T)$ 
  we obtain the following  weak formulation: 
\begin{align}
 &-\left(\beta\left(c^h_0\right),v^h\right)_{L^2(O^h;h)}\xi(0)-\int_0^T\left(\beta\left(c^h\right)(s),v^h\right)_{L^2(O^h;h)}\frac{d\xi}{dt}(s)\,ds\nonumber\\
 =&\int_0^T\left(c^h(s),\Delta^h v^h\right)_{L^2(O^h;h)}\xi(s)+\left(f\left(c^h,y^h\right)(s),v^h\right)_{L^2(O^h;h)}\xi(s)\,ds.\label{cdiscrweak}
\end{align} 
In particular, one can take $v^h:=v|_{\overline{O}^h}$ for  any $v\in C_0^{\infty}(U^h)$ as a test function in \eqref{cdiscrweak}. 
Thus, using the interpolation operators and properties  \eqref{lift},    \eqref{PiProp}, and \eqref{Lambda0} (\eqref{nLambda0})  one deduces from the space-discrete problem \eqref{PMdiscr}-\eqref{SDEdiscr} an approximation to the original space-continuous system \eqref{PM}-\eqref{SDE}:  
\begin{subequations}\label{defd}
\begin{empheq}[left={\empheqlbrace\,}]{align}
&-\left(\beta\left(\Pi^hP^hc_0\right),\Pi^hv\right)_{L^2(O)}\xi(0)-\int_0^T\left(\beta\left(\Pi^hc^h\right)(s),\Pi^hv\right)_{L^2(O)}\frac{d\xi}{dt}(s)\,ds\nonumber\\
 &=\int_0^T\left(\Pi^hc^h(s),\Pi^h\Delta^h v\right)_{L^2(O)}\xi(s)+\left(f\left(\Pi^hc^h,\Pi^hy^h\right)(s),\Pi^hv\right)_{L^2(O)}\xi(s)\,ds,\nonumber\\&\text{for all }v\in C_0^{\infty}(U^h),\ \xi\in C_0^{\infty}[0,T)\quad\text{ a.s. in }\Omega,\label{weakPch}\\
 & \Lambda^hc^h=0\quad\text{ in } (0,T)\times \partial O\qquad (\partial_{\nu} \Lambda^hc^h=0\quad\text{ a.e. in }(0,T)\times\partial O)\quad \text{ a.s. in }\Omega,\label{bcappr}\\
 &\Pi^h y^h(s)=\Pi^hP^hy_0+\int_0^ta\left(\Pi^h y^h\right)(s)\,d W(s)+\int_0^tb\left(\Pi^hc^h,\Pi^h y^h\right)(s)\,ds\nonumber\\
 &\text{in }L^2(O)\text{ for all }t\in[0,T]\quad\text{ a.s. in }\Omega.\label{Pyh}
\end{empheq}
\end{subequations}
This system will be analysed in the subsequent subsections.

\subsection{A priori estimates}
It is easy to see that for  solutions of the semi-discrete system \eqref{PMdiscr}-\eqref{SDEdiscr} one has a set of uniform in $h$ estimates which is very similar to estimates \eqref{estc1}-\eqref{esty1} for the original space-continuous system  \eqref{PM}-\eqref{SDE}. Indeed, one just needs to replace the spatial differential operators by their discrete versions  \eqref{dkh}-\eqref{nuh}, integration over $O$ by summation over $O^h$, and the $L^p(O)$ norms by their discrete analogies $L^p(O^h;h)$  and to use \eqref{EstPr1}-\eqref{EstPr2} in order to estimate $\left(c^h_0,y^h_0\right)$, the  discrete version of the Green's first identity \eqref{GrI}, 
the boundary conditions \eqref{bcdiscr}, as well as relation \eqref{chain} 
instead of the chain rule while dealing with spatial derivatives. Moreover, since we assumed that $\beta\in C^1(\R_0^+)$, the uniform  estimates for $c_{h}$ can be directly transformed into the corresponding uniform estimates for $\beta(c_{h})$. This spares the need of constructing  more  complicated  transformations such as those derived in {\it Subsections \ref{secPMTr}-\ref{SecPsi}}. We thus get the following sets of estimates:
\begin{subequations}\label{estch}
 \begin{align}
 &\left\|\beta\left(c^h\right)\right\|_{L^{\infty}(\Omega;L^{\infty}((0,T);L^{\infty}(\overline{O}^h;h)))}\leq \Cl{Ch},
 \\
  &\left\| \beta\left(c^h\right)\right\|_{L^{\infty}(\Omega;L^{\infty}((0,T);H^1(\overline{O}^h;h)))}\leq \Cr{Ch},\label{nabbetah}
 \\
 &\left\|\partial_t \beta\left(c^h\right)\right\|_{L^{\infty}(\Omega;L^2((0,T);L^2(\overline{O}^h;h)))}\leq \Cr{Ch},
 \\
 &\left\|D_r \beta\left(c^h\right)\right\|_{L^{2}(r,T;L^2(\Omega;L^2(\overline{O}^h;h)))}\leq \Cr{Ch},
 \\
 &\left\|h_2^{-\frac{1}{2}}\left(D_{r+h_2}\beta\left(c^h\right)-D_r \beta\left(c^h\right)\right)\right\|_{L^2((r,T-h_2)\times \Omega;L^2(\overline{O}^h;h))}
 \leq \Cr{Ch}, 
 \\
 &\left\|\partial_t D_r \beta\left(c^h\right)\right\|_{L^2((r,T)\times \Omega;H^{-2}(\overline{O}^h;h))}\leq\Cr{Ch}
\end{align}
and
\end{subequations}
\begin{subequations}\label{estyh}
\begin{align}
&\left\|y^h\right\|_{C^{\frac{1}{2}}([0,T];L^q(\Omega;L^q(\overline{O}^h;h)))}\leq\Cr{Ch},\\
&\left\|y^h\right\|_{L^{\infty}((0,T);L^2(\Omega;H^1(\overline{O}^h;h)))}\leq \Cr{Ch},\label{nabyh}\\
&\left\|D_ry^h\right\|_{L^{\infty}(r,T;L^2(\Omega;L^2(\overline{O}^h;h)))}\leq\Cr{Ch},\\
&\left\|h_2^{-\frac{1}{2}}(D_{r+h_2}y^h-D_r y^h)\right\|_{L^{\infty}(r,T-h_2;L^2(\Omega;L^2(\overline{O}^h;h)))}\leq\Cr{Ch}\\
&\left\|h_1^{-\frac{1}{2}}(D_r y^h(\cdot+h_1)-D_r y^h)\right\|_{L^{\infty}(r,T-h_1;L^2(\Omega;L^2(\overline{O}^h;h)))}
 \leq\Cr{Ch}.
\end{align}
\end{subequations}
Combining \eqref{estch}-\eqref{estyh} with \eqref{TransL1}-\eqref{TransL3} we conclude that the pair $\left(\Lambda^h\beta\left(c^h\right),\Lambda^h\left(y^h\right)\right)$ satisfies  
\begin{subequations}\label{estchL}
 \begin{align}
 &\left\|\Lambda^h\beta\left(c^h\right)\right\|_{L^{\infty}(\Omega;L^{\infty}((0,T);L^{\infty}(O)))}\leq \Cl{ChL},
 \\
  &\left\|\Lambda^h\beta\left(c^h\right)\right\|_{L^{\infty}(\Omega;L^{\infty}((0,T);H^1(O)))}\leq \Cr{ChL},\label{nabLbh}
 \\
 &\left\|\partial_t \Lambda^h\beta\left(c^h\right)\right\|_{L^{\infty}(\Omega;L^2((0,T);L^2(O)))}\leq \Cr{ChL},
 \\
 &\left\|D_r \Lambda^h\beta\left(c^h\right)\right\|_{L^{2}(r,T;L^2(\Omega;L^2(O)))}\leq \Cr{ChL},
 \\
 &\left\|h_2^{-\frac{1}{2}}\left(D_{r+h_2}\Lambda^h\beta\left(c^h\right)-D_r \Lambda^h\beta\left(c^h\right)\right)\right\|_{L^2((r,T-h_2)\times \Omega;L^2(O))}
 \leq \Cr{ChL}, 
 \\
 &\left\|\partial_t D_r \Lambda^h\beta\left(c^h\right)\right\|_{L^2((r,T)\times \Omega;H^{-2}(U^h))}\leq\Cr{ChL},
\end{align}
\end{subequations}
\begin{subequations}\label{estyhL}
\begin{align}
&\left\|\Lambda^h y^h\right\|_{C^{\frac{1}{2}}([0,T];L^q(\Omega;L^q(O)))}\leq\Cr{ChL},\\
&\left\|\Lambda^h y^h\right\|_{L^{\infty}((0,T);L^2(\Omega;H^1(O)))}\leq \Cr{ChL},\label{nabLyh}\\
&\left\|D_r\Lambda^h y^h\right\|_{L^{\infty}(r,T;L^2(\Omega;L^2(O)))}\leq\Cr{ChL},\\
&\left\|h_2^{-\frac{1}{2}}(D_{r+h_2}\Lambda^h y^h-D_r \Lambda^h y^h)\right\|_{L^{\infty}(r,T-h_2;L^2(\Omega;L^2(O)))}\leq\Cr{ChL}\\
&\left\|h_1^{-\frac{1}{2}}(D_r \Lambda^h y^h(\cdot+h_1)-D_r \Lambda^h y^h)\right\|_{L^{\infty}(r,T-h_1;L^2(\Omega;L^2(O)))}
 \leq\Cr{ChL}.
\end{align}
\end{subequations}
Further, \eqref{nabbetah} and  \eqref{nabyh} together with \eqref{compare} and   assumption $\beta^{-1}\in C^1(\R_0^+)$ (compare \eqref{betamm}), yield that
\begin{align}
&\left\|\Pi^hc^h-\Lambda^h c^h\right\|_{L^{\infty}(\Omega;L^{\infty}((0,T);L^2(O)))}\leq h\Cl{Cdif1},\label{PiLch}\\
 &\left\|\Pi^h\beta\left(c^h\right)-\Lambda^h \beta\left(c^h\right)\right\|_{L^{\infty}(\Omega;L^{\infty}((0,T);L^2(O)))}\leq h\Cr{Cdif1}\label{PiLbeh}
\end{align}
and 
\begin{align}
 \left\|\Pi^hy^h-\Lambda^h y^h \right\|_{L^{\infty}((0,T);L^2(\Omega;L^2(O)))}\leq h\Cr{Cdif1}.\label{PiLyh}
\end{align}
Finally, due to \eqref{estPiP} we also have for the approximations of the initial data the estimates
\begin{align}
 &\left\|\Pi^hP^h c^h_0-c_0\right\|_{L^2(O)}\leq h\Cl{C14},\label{estc0h}\\
 &\left\|\Pi^hP^h y^h_0-y_0\right\|_{L^2(O)}\leq h\Cr{C14}.\label{esty0h}
\end{align}
\subsection{Compactness and convergence}
Thanks to estimates \eqref{estchL}-\eqref{estyhL} we are once again in the position when we can directly apply \cite[Theorem 2]{BallySaussereau}. It  yields that
\begin{align}
 \left\{\left(\Lambda^h\beta\left(c^h\right),\Lambda^h y^h\right):\quad h=\frac{1}{M+1},\ M\in\N\right\}\quad\text{ is precompact in }\left(L^2((0,T)\times O\times\Omega)\right)^2.\label{LamdaComp}
\end{align}
 Together with estimates \eqref{PiLbeh}-\eqref{PiLyh} this leads to 
\begin{align}
 \left\{\left(\Pi^h\beta\left(c^h\right),\Pi^h y^h\right):\quad h=\frac{1}{M+1},\ M\in\N\right\}\quad\text{ is precompact in }\left(L^2((0,T)\times O\times\Omega)\right)^2.\label{PiComp1}
\end{align}
Further, using $\beta^{-1}\in C^1(\R_0^+)$ and property \eqref{PiProp} for $\varphi:=\beta$, we conclude with \eqref{PiComp1} that
\begin{align}
 \left\{\left(\Pi^h c^h,\Pi^h y^h\right):\quad h=\frac{1}{M+1},\ M\in\N\right\}\quad\text{ is precompact in }\left(L^2((0,T)\times O\times\Omega)\right)^2.\label{PiComp}
\end{align}
Combining estimates \eqref{estchL}-\eqref{estyhL}, and \eqref{PiLch}, \eqref{PiLyh} with the compactness results \eqref{LamdaComp} and \eqref{PiComp} and using  assumption $\beta^{-1}\in C^1(\R_0^+)$  we deduce that  there exists a  pair of functions $(c,y):[0,T]\times \overline{O}\times\Omega\rightarrow\R^+_0\times \R^+_0$ which satisfies conditions \ref{Defsol1}-\ref{Defsol5} from {\it Definition \ref{Defsol}}   
and a sequence $h_n\underset{n\rightarrow\infty}{\rightarrow}0$ such that 
\begin{align}
 &\left(\Pi^{h_n}c^{h_n},\Pi^{h_n}y^{h_n}\right)\underset{n\rightarrow\infty}{\rightarrow}(c,y)\quad\text{ in }L^2((0,T)\times O\times\Omega),\quad\text{ in }L^2((0,T)\times O)\text{ a.s. in }\Omega\label{convcy}
 \end{align}
 and
\begin{align}
 &\Lambda^{h_n}c^{h_n}\underset{n\rightarrow\infty}{\overset{*}{\rightharpoonup}}c\quad\text{ in }L^{\infty}(\Omega;L^{\infty}((0,T);H^1(O))).\label{convnch}
\end{align}
Due to \eqref{estc0h}-\eqref{esty0h} it holds also  that
\begin{align}
 \left(\Pi^{h_n}P^{h_n}c^{h_n}_0,\Pi^{h_n}P^{h_n}y^{h_n}_0\right)\underset{n\rightarrow\infty}{\rightarrow}(c_0,y_0)\quad\text{ in }(L^2(O))^2.\label{convc0y0}
\end{align}
Using \eqref{convcy} and \eqref{convc0y0}, assumption $\beta,f,a,b\in C^1$, and the fact that
\begin{align}
 \Pi^{h_n}\Delta^{h_n}v\underset{n\rightarrow\infty}{\rightarrow} \Delta v\quad \text{ in }L^2(O)\nonumber
\end{align}
due to \eqref{PiDelh}, we can pass to the limit in the weak formulation \eqref{weakPch} and the SDE \eqref{Pyh}. 
Finally, combining \eqref{convnch} with the boundary conditions \eqref{bcappr} and  the trace theorem, we conclude that $c$ satisfies either of the boundary conditions. Altogether, this means that $(c,v)$ is a weak solution to the PDE \eqref{PM} and a strong solution to the SDE \eqref{SDE}, and it has  the desired  regularity. 
Thus, we proved the following Lemma:
\begin{Lemma}[Existence for a nondegenerate case]\label{ExistLemma}
Let  {\it Assumptions \ref{Assump1}} be satisfied. Assume in addition that $\beta\in C^2(\R_0^+)$ Then there exists  a weak-strong global solution in terms of {\it Definition \ref{Defsol}}  to system \eqref{PM}-\eqref{SDE}.
\end{Lemma}

\section{Proof of the existence {\it Theorem \ref{ExistTh}}}\label{existence}
In this final section we discuss how our compactness {\it Theorem \ref{CompTh}} can be used in order to prove the existence of solutions to the original degenerate system \eqref{PM}-\eqref{SDE}. Let $T>0$ be arbitrary. Choose $R_0,R_1>0$ large enough so that $(\beta,f,a,b,c_0,y_0)\in{\mathcal{P}}(T,(R_0/2,R_1/2))$. For each $\varepsilon\in(0,1)$ let $\beta_{\varepsilon}$ be such that
\begin{align}
&(\beta_{\varepsilon},f,a,b,c_0,y_0)\in{\mathcal{P}}(T,(R_0,R_1)),\label{class}\\
&\beta_{\varepsilon}\in C^2(\R_0^+),\nonumber\\
&\beta_{\varepsilon}\underset{\varepsilon\rightarrow0}{\rightarrow}\beta\quad\text{ in }C[0,R_2],\nonumber
\end{align}
where $R_2$ is an upper bound for the $L^{\infty}$-norm of $c$, compare \eqref{R2}.
We introduce a family of approximating  problems
 \begin{subequations}\label{epsPME}
 \begin{empheq}[left={\text{in } \Omega\ \empheqlbrace\,}]{align}
&\partial_t \beta_{\varepsilon}(c_{\varepsilon})=\Delta c_{\varepsilon}+f(c_{\varepsilon},y_{\varepsilon})&&\text{ in }(0,T]\times O,\label{epsc}\\
 &c_{\varepsilon}=0\qquad (\partial_{\nu}c_{\varepsilon}=0)&&\text{ in }(0,T]\times \partial O,\label{epsbc}\\
 &c_{\varepsilon}=c_0&&\text{ in }\{0\}\times O,\label{epsinic0}
\end{empheq}
\end{subequations}
\begin{subequations}\label{epsSDE}
 \begin{empheq}[left={\text{in } O\ \empheqlbrace\,}]{align}
 &dy_{\varepsilon}=a(y_{\varepsilon})\,d W+b(c_{\varepsilon},y_{\varepsilon})\,dt &&\text{ in } (0,T]\times\Omega,\label{epsy}\\
 &y_{\varepsilon}=y_0&&\text{ in } \{0\}\times\Omega.\label{epsiniy0}
\end{empheq}
\end{subequations}
 The existence of  a solution $\left(c_{\varepsilon},y_{\varepsilon}\right)$ in terms of {\it Definition \ref{Defsol}} to system \eqref{epsPME}-\eqref{epsSDE}  follows from {\it Lemma \ref{ExistLemma}} of the previous section. Moreover, due to assumption \eqref{class} and the compactness {\it Theorem \ref{CompTh}}  family $\left\{\left(c_{\varepsilon},y_{\varepsilon}\right)\right\}_{\varepsilon\in(0,1)}$ is precompact in $\left(L^2((0,T)\times O\times\Omega)\right)^2$. Standard compactness arguments and a limit procedure  then yield 
the existence of a solution in terms of {\it Definition \ref{Defsol}} to the original degenerate system \eqref{PM}-\eqref{SDE} for any $T>0$. 

Finally, we observe that the particular choice of the zero starting time was not essential for our arguments heretofore.
Indeed, all previous results continue to hold if we replace the interval $[0,T]$ by $[t_0,t_0+T]$ for any $t_0>0$  and consider the 'shifted' in time filtration  $({\mathcal F_{t+t_0}})_{t\geq0}$ instead of the original one. Thus, we  obtain a solution defined for all times by defining it successively in $[0,T]$, $[T,2T]$, and so on. The proof of {\it Theorem \ref{ExistTh}} is thus complete. 
\begin{Remark}[Numerical schema] As a by-product of the   constructions from this Section and the  previous {\it Section \ref{SecSemi}}  we have a numerical schema for  \eqref{PM}-\eqref{SDE} which is, at least theoretically, converging.
\end{Remark}

\phantomsection
\addcontentsline{toc}{section}{References}
\printbibliography

\end{document}